\documentclass{amsart}
\usepackage{ifthen}
\usepackage{hyperref}

\usepackage{amsmath,amssymb,amsthm,verbatim,paralist,graphicx}

\usepackage[active]{srcltx}
\usepackage[T1]{fontenc}
\usepackage[latin1]{inputenc}
\usepackage{yfonts}
\usepackage{color}
\DeclareMathOperator{\rk}{r}
\DeclareMathOperator{\ind}{ind}
\DeclareMathOperator{\cl}{cl}
\DeclareMathOperator{\supp}{supp}
\DeclareMathOperator{\sep}{sep}
\DeclareMathOperator{\Dir}{Dir}

\numberwithin{equation}{section}
\newtheorem{theorem}{Theorem}[section]
\newtheorem{prop}{Proposition}[section]
\newtheorem{lemma}{Lemma}[section]
\newtheorem{corollary}{Corollary}[section]
\newtheorem{definition}{Definition}[section]

\newtheorem{remark}[lemma]{Remark}
\newtheorem{problem}{Problem}[section]

\newtheoremstyle{TheoremNum}
        {\topsep}{\topsep}              
        {\itshape}                      
        {}                              
        {\bfseries}                     
        {.}                             
        { }                             
        {\thmname{#1}\thmnote{ \bfseries #3}}
    \theoremstyle{TheoremNum}
    \newtheorem{thmn}{Theorem}

\begin{document}
\title{Lexicographic Extensions preserve Euclideaness}  \author{Winfried Hochst\"attler} \author{Michael Wilhelmi}
\address{FernUniversit\"at in Hagen}
\email{Winfried.Hochstaettler@FernUniversitaet-Hagen.de}
\email{mail@michaelwilhelmi.de}

\begin{abstract}
We prove that a lexicographical extension of a Euclidean oriented
  matroid remains Euclidean. 
\end{abstract}

\dedicatory{}

\maketitle


\section{Introduction}
Linear Programming and in particular the simplex algorithm can be generalized to the combinatorial setting of oriented matroids \cite{2}. Robert Bland, the father of the classical pivoting rule that prevents cycling in linear programs with degenerate vertices, is also one of the fathers of oriented matroids. In \cite{Blandpivot} he already remarks that ``Most proofs of finiteness (of the simplex method) ... invoke monotonicity of the function value. This property cannot, as far as we know, be nicely translated to oriented matroids.''
A little later Edmonds and Fukuda~\cite{6} presented an example showing that non-degenerate cycling indeed is possible in non-linear oriented matroid programs. Edmonds and Mandel~\cite{5} showed that the reason for this phenomenon lies in the non-existence of certain parallel extensions. Therefore, they named  the oriented matroids, where non-degenerate cycling cannot occur {\em Euclidean}.

Euclideanness of oriented matroid programs is far from well understood.
While it is clear that a minor of a Euclidean oriented matroid program is still Euclidean, see \cite{5}, 9.III.17 or \cite{2}, 10.5.6, it is an open question which kind of single-element extensions of oriented matroids preserve Euclideanness. Our interest in Euclidean oriented matroid programs has its origin in the fact that they provide the possibility to sweep the feasible region with a hyperplane. Intuitively, it is clear that this should lead to a vertex-shelling of the feasible region. We will show in our companion paper~\cite{Paper2} that this is indeed possible. In order to make the recursion in the shelling procedure work, though, we were confronted with the problem that the cone over a vertex in an oriented matroid is in general no longer a bounded region in the minor. In order to fix this, we had to consider extensions of the original program.

While it was known that a parallel extension of a Euclidean oriented matroid program preserves Euclideanness, for other extensions this is an open question. It is clear that parallel extensions will not suffice to bound the region in the contracted matroid program. Here, we will show that lexicographic extensions preserve Euclideanness. 
More precisely, we will prove:

\begin{theorem}\label{cor:lexExtStaysEucl2a}
Let $\mathcal{O}$ be a Euclidean oriented matroid and let $\mathcal{O} \cup p$ be a lexicographic extension of $\mathcal{O}$.
Then $(\mathcal{O} \cup p,p,f)$ is a Euclidean matroid program.
\end{theorem}

In Section~\ref{sec:two} of this paper we present some preliminaries,
required definitions, the proof that a lexicographically extended
oriented matroid program is bounded and that Euclideaness is preserved
if we exchange the target function and the hyperplane at infinity.
Section~\ref{sec:three} is then devoted to the proof of Theorem
\ref{cor:lexExtStaysEucl2a}. We actually verify its contraposition.
We compare the directions of edges in the extended oriented matroid
program with those in some programs of the original oriented matroid
and prove that, if there is a directed cycle in the extended matroid,
we always obtain a directed cycle in an oriented matroid program of
the original matroid.  
Based on that result, we also show the general case.

\begin{theorem}\label{theorem:SecondMainTheorem}
Any lexicographic extension of a Euclidean oriented matroid
is Euclidean.
\end{theorem}

We prove this theorem in Section~\ref{sec:four}. Again, a directed
cycle in an extended matroid yields a directed cycle in the old one,
but here we have to find that old cycle. The difficulty that might
arise is when the directed cycle uses an edge that does not correspond
to an edge in the original oriented matroid. In order to replace that edge, we convert the new
cycle into a simpler form, a so-called {\em normalized cycle} from
which we can derive either a so-called {\em corresponding directed
  cycle} in the original matroid or, with the help of a so-called {\em
  Projection Lemma}, we can project it to a directed cycle at the
plane of infinity, reducing the proof to Theorem
\ref{cor:lexExtStaysEucl2a}.  

\section{Premilinaries}\label{sec:two}

For the basic concepts/notations of oriented matroid programming, we
refer to \cite{2}, Chapter 10. Unlike in \cite{2}, we call an oriented
matroid program $(\mathcal{O},g,f)$ {\em bounded} if $X_g = +$ for
every cocircuit $X$ of the feasible region. 
We
use the notion of Euclideaness of an oriented matroid program in
\cite{2}, 10.5.2. An oriented matroid program $(\mathcal{O},g,f)$ is Euclidean if and only if its directed cocircuit graph $G_f$ contains no directed cycles.  Additionally, we call an
oriented matroid $\mathcal{O}$ {\em Euclidean} if all programs
$(\mathcal{O},g,f)$ with $f \neq g$ of the groundset ($f$ not a
coloop, $g$ not a loop) are Euclidean and we call it {\em Mandel} if
there exists a single element extension $\mathcal{O}'
= (\mathcal{O} \cup f)$ of $\mathcal{O}$ in general position such that
$(\mathcal{O}',g,f)$ is Euclidean for all elements $g$ of the
groundset that are not loops.  

\subsection{Directing Edges in the Graph $G_f$ and Beyond}

Recall, that if $X$ and $Y$ are covectors in an oriented matroid
$\mathcal{O}$ we have $z(X \circ Y) = z(Y \circ X) = z((-X) \circ Y)$
for their zero sets. Moreover, if $X,Y$ are cocircuits with $X \neq
\pm Y$ then $\supp(Y) \setminus \supp(X) \ne \emptyset \ne \supp(X)
\setminus \supp(Y)$. A covector $F$ is an {\em edge} in $\mathcal{O}$
if $z(F)$ is a coline in the underlying matroid
$\mathcal{M}(\mathcal{O})$.  We say a cocircuit $X$ is {\em lying on the edge $F$} or {\em on the line spanned by} (or {\em through}) $F$ if $z(X) \supset z(F)$. We call two covectors a {\em modular
  pair} or {\em comodular} in $\mathcal{O}$ if their zero-sets form a
modular pair in $\mathcal{M}(\mathcal{O})$.  Hence two cocircuits $X$
and $Y$ are comodular if and only if $X \circ Y$ is an edge in~$\mathcal{O}$.  Also we will call $(X,Y)$ an {\em
  edge} in $\mathcal{O}$ if $X$ and $Y$ are conformal comodular cocircuits.

The following propositions are immediate.

\begin{prop}\label{prop:CocircuitsLyingOnEdges1}
Let $F$ be an edge in an oriented matroid $\mathcal{O}$ and
$X \neq \pm Y$ cocircuits and  $z(F) \subseteq z(X)$ as well as  $z(F) \subseteq z(Y)$. Then
$z(X \circ Y) = z(F).$
\end{prop}

\begin{prop}\label{prop:unqieCocircuitOnEndge}
  Let $F$ be an edge of an oriented matroid $\mathcal{O}$ and $e \in
  \supp(F)$.  Then there is, up to sign reversal, a unique cocircuit
  $Z$ with $e \not \in \supp(Z) \subseteq \supp(F)$.
\end{prop}

\begin{prop}[Uniqueness Of Cocircuit Elimination]\label{prop:UniqueCoEl}
  Let $X,Y$ be comodular cocircuits of an oriented matroid and $e \in
  \sep(X,Y)$. Then cocircuit elimination of $e$ between $X$ and $Y$ yields a unique cocircuit $Z$ such that $\supp(Z) \subseteq \supp(X
  \circ Y) \setminus e$. Furthermore, we have $Z_f=X_f\circ Y _f$ for
  all $f \notin \sep(X,Y)$.
\end{prop}

We introduce a direction for modular pairs of
cocircuits.

\begin{definition}\label{def:Dirfunction}
Let $(\mathcal{O},g,f)$ be an oriented matroid program.
Let 
\[ \mathcal{P} = \{ (X,Y) \mid X \text{ and } Y \text{ are comodular cocircuits in } \mathcal{O} \text{ with } X_g = Y_g \neq 0 \}.\]

We say $(X,Y) \in \mathcal{P}$ is directed from $X$ to $Y$ and write
$X \to_{g,f} Y$ \\ (or $X \to_{(\mathcal{O},g,f)} Y$), if cocircuit
elimination of $g$ between $-X$ and $Y$ yields $Z$ satisfying $Z_f=+$.
We write $X \leftarrow_{g,f} Y$ if $Z_f=-$ and $X
\leftrightarrow_{g,f} Y$ if $Z_f=0$ (then we say $(X,Y)$ is
\em{undirected}).  We may omit the subscript, in case it is obvious.

\end{definition} 

The following again is immediate (we use the order $- \prec 0 \prec +$ in $\{-,0,+\}$)::

\begin{prop} \label{prop:basicDirectionProps} Let $(X,Y) \in \mathcal{P}$. Then it holds
  \begin{enumerate} 
\item $X \to Y \iff -Y \to -X$.
   \item  $X \to_{(\mathcal{O},g,f)} Y \Leftrightarrow X \to_{({}_{-g}\mathcal{O},g,f)} Y$ and $X \to_{(\mathcal{O},g,f)} Y \Leftrightarrow Y \to_{({}_{-f}\mathcal{O},g,f)} X$ where ${}_{-g}\mathcal{O}$ resp.  ${}_{-f}\mathcal{O}$ are derived from $\mathcal{O}$ reorienting $g$ resp. $f$.
    \item If $X_f \prec Y_f$, then $X \to
    Y$. If $X_f = Y_f = 0$, then $X \leftrightarrow Y$.

  \item Let $e \in \sep(X,Y)$, $X \to Y$
    and cocircuit elimination of $e$ between $X$ and $Y$ yield $W$.
    Then $(X,W),(W,Y) \in \mathcal{P}$ and 
    $X \to W$ holds as well as $W \to Y$.
   \item If $(X,Y)$ is undirected, then $(Y,X), (-X,-Y)$ and $(-Y,-X)$ are also undirected in every reorientation of $\mathcal{O}$.
  \end{enumerate}
\end{prop}

We add a small Proposition.
\begin{prop} \label{prop:TransitivityOfDir}
Let $X_g = Y_g \neq 0$ and $X_f = Y_f \neq 0$ and $h \in E \setminus \{f,g\}$. If $X \leftrightarrow_{(\mathcal{O},g,f)} Y$, then $X \to_{(\mathcal{O},g,h)} Y$ implies $X \to_{(\mathcal{O},f,h)} Y$.
\end{prop} 
\begin{proof}
Cocircuit elimination of $g$ between $-X$ and $Y$ yields $Z$ with $Z_f = 0$. Hence, elimination of $f$ between $-X$ and $Y$ yields again $Z$.
\end{proof}

If $X$ and $Y$ are a conformal comodular pair, the former
definition is compatible with the orientation of the cocircuit graph
of an oriented matroid program.

\begin{definition} [see \cite{2}, 10.1.16]\label{def:Graph}
Let $(\mathcal{O},g,f)$ be an oriented matroid program. 
Let $G_f$ be the graph whose vertices are the cocircuits of $(\mathcal{O},g)$ and whose edges are the edges $(Y^1,Y^2)$ between two comodular, conformal cocircuits.
The edges in $G_f$ are (un-)directed like in $(\mathcal{O},g,f)$.
We call the graph $(G \setminus f)_f$ the directed subgraph of $G_f$
induced by the cocircuits of $(\mathcal{O} \setminus f,g)$.  If two
comodular cocircuits $Y^1,Y^2$ with $Y^1 \circ Y^2$ are conformal in
$(\mathcal{O} \setminus f,g)$ but not in $G_f$, then we have $Y^1_f \succ
Y^2_f$ (or $Y^1_f \prec Y^2_f$) in $\mathcal{O}$ and we direct the edge 
$Y^1 \rightarrow Y^2$ (or $Y^1 \leftarrow Y^2$) in $(G \setminus f)_f$.
\end{definition}

We define paths in $G_f$. We give a reformulation of Definition 10.5.4 in \cite{2}.

\begin{definition}
A {\em path} in $G_f$ from $P^1$ to $P^n$ is a sequence
\[ P = [P^1, \hdots, P^n] \]
of vertices of $G_f$ such that $(P^i,P^{i+1})$ is an edge in $G_f$ for all $1 \leq i < n$.
A path $P$ is {\em undirected} if $P^i \leftrightarrow P^{i+1}$ for
all $1 \leq i < n$ and {\em directed} if is not undirected and either
$P^i \to P^{i+1}$ or  $P^i \leftrightarrow P^{i+1}$ for all $1
\leq i < n$, in which case we say that the path {\em is
  directed from $P^1$ to $P^n$}.

A path is {\em closed} if $P^1 = P^n$. A closed path is a {\em cycle}
if $P^1, \hdots, P^{n-1}$ are pairwise distinct.
\end{definition}

We collect some more well-known facts.

\begin{prop}  \label{prop:PathGraphWellKnownFacts}
  \begin{enumerate}
  \item \label{prop:ClosedPathContainsCycle}
Let $P = [P^1, \hdots, P^n]$ be a closed path in a graph and and let $X,Y$ be subsequent vertices of $P$. 
Then there always exists a subpath $Q$ of $P$ containing $X$ and $Y$ that is a cycle.
\item \label{pathsOnEdges} If $X,Y$ are comodular cocircuits in $G_f$
  then there is a path $X^0, \hdots, X^n$,  $X=X_0, Y=X_n$ in $G_f$ with
  distinct cocircuits such that $\supp(X_i \circ X_j) = \supp(X \circ
  Y)$ and all pairs $(X_i, X_j)$ are directed as $(X,Y)$ for all $i < j$.
\item \label{prop:DirectedPaths}
Let $P$ be a directed path $P = (X,X^1, \hdots, X^n, Y)$ between two cocircuits $X$ and $Y$ in $G_f$.
If $X_f \prec Y_f$ the path 
is directed
from $X$ to $Y$.
\item \label{prop:PathGraphWellKnownFacts_directedCycle} Let $(\mathcal{O},g,f)$ be an
  oriented matroid program that contains a directed cycle in the graph
  $G_f$. Then all cocircuits in the cycle have the same $f$-value, in particular $X_f \ne 0$ holds everywhere.
  \end{enumerate}
\end{prop}

\begin{lemma}\label{prop:ChangeOfFAndG1}
  Let $(\mathcal{O},g,f)$ be an oriented matroid program and $X,Y$ be
  comodular cocircuits in $\mathcal{O}$ with $X_f = Y_f = X_g = Y_g
  \neq 0$.  Then 
  \[X \to_{g,f} Y \Leftrightarrow X \leftarrow_{f,g} Y,\] i.e.\ the
  orientation switches if we interchange the roles of the objective
  function and infinity.
\end{lemma}

\begin{proof} 
Possibly using Proposition \ref{prop:basicDirectionProps} (i), we may assume $X_f = +$.
Cocircuit elimination of $g$ between $-X$ and $Y$ yields $Z$ with $Z_g = 0$.
Cocircuit elimination of $f$ between $-X$ and $Y$ yields $W$ with $W_f = 0$.
We have three cases.
\begin{enumerate}
\item If $Z_f = 0$ then $Z = W$ hence $X \leftrightarrow Y$ in both
  programs.
\item If $Z_f = +$ let $h$ with $Y_h = 0$ and $X_h \neq 0$ we may
  assume $X_h = +$ hence $W_h = Z_h = -$.  Cocircuit elimination of
  $f$ between $-X$ and $Z$ yields $W^1$ with $W^1_g = -$ and $W^1_h =
  -$.  Since $X$ and $Y$ are comodular we must have $W = W^1$ hence
  $W_g = -$ and that is $X \leftarrow Y$ in  $G_g$.
\item If $Z_f = -$, we interchange the roles of $X$ and $Y$ and
(ii) yields $Y \leftarrow X$ in $G_g$.
\end{enumerate}
\end{proof}

Reorienting $g$, we derive a corollary:

\begin{corollary}\label{cor:ChangeOfFAndG2}
Let $(\mathcal{O},g,f)$ be an oriented matroid program and $X,Y$ be comodular cocircuits in $\mathcal{O}$ with $X_f = Y_f  = - X_g = - Y_g \neq 0$.
Then \[X \to_{g,f} Y \Leftrightarrow X \to_{f,g} Y.\] 
\end{corollary}

\begin{remark}

  The last lemma and its corollary immediately imply $X
  \leftrightarrow_{g,f} Y \Leftrightarrow X \leftrightarrow_{f,g} Y$.
\end{remark}

\begin{prop}\label{prop:DirectedCycleYieldsDirectedCycle}
If $\mathcal{G}$ is a directed cycle in $(\mathcal{O},g,f)$ then either $\mathcal{G}$ or $-\mathcal{G}$ is a directed cycle in $(\mathcal{O},f,g)$
\end{prop}

\begin{proof}
  All cocircuits $X$ of $\mathcal{G}$ satisfy $X_g=+$. By Proposition
  \ref{prop:PathGraphWellKnownFacts} (iv) all $X_f$ have the same value which is
  different from $0$. Hence, all edges in the circuit satisfy either
  the prerequisites of the Lemma~\ref{prop:ChangeOfFAndG1} or of its
  corollary. In the former case $\mathcal{G}$ is a directed cycle
  in $(\mathcal{O},f,g)$ as well, in the latter case  $-\mathcal{G}$ is (\see Proposition \ref{prop:basicDirectionProps}).
\end{proof}

We obtain immediately.

\begin{theorem} \label{theorem:EuclideanessStays}
If $(\mathcal{O},g,f)$ is a Euclidean oriented matroid program, so is $(\mathcal{O},f,g)$.
\end{theorem}

Clearly, reorienting elements from $E \setminus \{g,f\}$ does not
change the direction of any edge, and hence:

\begin{lemma}
Let $(\mathcal{O},g,f)$ be a Euclidean oriented matroid program with groundset $E$ and let
$S \subseteq E \setminus \{f,g\}$ and let  ${}_{-S}\mathcal{O}$ be derived from $\mathcal{O}$ by reorientation of the elements of $S$.
Then $({}_{-S}\mathcal{O},g,f)$ is a Euclidean oriented matroid program, too.
\end{lemma}
 
\subsection{Old and new cocircuits in single-element extensions} 

For the notion of single-element extensions of oriented matroids and its corresponding localizations, we refer to \cite{2}, Chapter 7.1.
We call an extension non-trivial if it is not an extension with a coloop. We mention one fact explicitly.

\begin{theorem}[see \cite{2}, 7.1.4]\label{theorem:ext}
Let $\mathcal{O}$ be an oriented matroid of rank $\rk$ and $\mathcal{C}$ the set of its cocircuits,
let $\mathcal{O}' = \mathcal{O} \cup p$ be a non-trivial 
single-element extension of $\mathcal{O}$
and let $\sigma$ be the corresponding localization. Let 
\[  S =  \{Z = (Y,\sigma(Y)) \colon Y \in \mathcal{C} \}.  \] 
Then the cocircuits $ \mathcal{C}'$ of $\mathcal{O}'$ are given by
\begin{align*} \mathcal{C}' = S \cup \{ Z' = (Y^1 \circ Y^2,0) \colon & Y^1, Y^2 \in \mathcal{C}, \sigma(Y^1) = -\sigma(Y^2) \neq 0,  \\
 & \sep(Y^1,Y^2) = \emptyset, \rk(z(Y^1 \circ Y^2)) = \rk -2 \} 
 \end{align*}
 We say that the cocircuit $Z \in S$ is {\em derived from} the
 cocircuit $Y$ or {\em old} and the cocircuit $Z' \in \mathcal{C}'
 \setminus S$ is {\em derived from the edge} $Y^1 \circ Y^2$ or {\em
   from the cocircuits} $Y^1,Y^2$ or {\em new}. We call an edge of two
 old cocircuits old as well. We call an edge new if at least one of its cocircuits is new.
 \end{theorem}
 Let us mention some well-known relations between edges and vertices
 in the oriented matroid and its extension.
 \begin{prop} \label{prop:uniqueEdges2} Let $\mathcal{O}$ be an
   oriented matroid and let $\mathcal{O}' = \mathcal{O} \cup p$ be a
 non-trivial single-element extension of $\mathcal{O}$.
  \begin{enumerate}
  \item \label{prop:uniqueEdges} 
    If $Z = (X \circ Y,0)$ is a cocircuit in $\mathcal{O}'$
    derived from the edge $X \circ Y$ where $X,Y$ are conformal
    comodular cocircuits in $\mathcal{O}$, then the pair $\{X,Y\}$ is
    unique.
  \item\label{prop:cocircuitsInSingleElementExtensions} If $Y$ is a
    cocircuit of $\mathcal{O}'$ with $Y_p \neq 0$, then $Y \setminus
    p$ is a cocircuit in $\mathcal{O}$, $Y$ is old.
\item If $p$ is in general position, then every cocircuit $X$ of $\mathcal{O}'$ with $X_p = 0$ is new.
  \end{enumerate}
\end{prop}

\subsection{Lexicographic extensions}

We start with a technical definition.

\begin{definition}
  Let $I = [e_1, \hdots, e_k ]$ be an ordered subset of the groundset
  of an oriented matroid $\mathcal{O}$ and let $Y$ be a cocircuit in
  $\mathcal{O}$. We call the index $i$ such that $Y_{e_i} \neq 0$ and
  $Y_{e_j} = 0$ for all $j < i$ the {\em index of the cocircuit $Y$
    wrt.\ I}. If $Y_{e_i} = 0$ for all $1 \leq i \leq
  k$, let the index of the cocircuit wrt.\ I be $k+1$.
\end{definition}

We define a lexicographic extension as follows.

\begin{definition}[see \cite{2}, 7.2.4] \label{prop:lexExtBjoerner}
Let $\mathcal{O}$ be an oriented matroid with  $\mathcal{C}$ being its set of cocircuits, $I = [e_1, \hdots, e_k]$ an ordered subset of the groundset, and $\alpha = [\alpha_1, \hdots, \alpha_k] \in \{+,-\}^k.$
Then the lexicographic extension $\mathcal{O}[I^{\alpha}] = \mathcal{O}[e^{\alpha_1}_1, \hdots, e^{\alpha_k}_k]$ of $\mathcal{O}$ is given by the localization
$\sigma \colon \mathcal{C} \rightarrow \{+,-,0\}$ with
\[ \sigma(Y) = \begin{cases}\alpha_iY_{e_i} & \text{if for the index } i \text{ of }Y \text{ wrt.\ } I \text{ holds } i \leq k , \\ 0 &  otherwise. \end{cases}\]
We call the lexicographic extension \em{positive} and write $\mathcal{O}[I^+]$ iff $\alpha_i = +  $ for all $i \leq k$.
\end{definition}

It is immediate that a lexicographic extension can always be written
as a reorientation of the original, a positive lexicographic extension
and then reversing the reorientation in the extension.  Until the end
of this chapter, let $\mathcal{O}$ always be an oriented matroid of
rank $\rk \geq k$, let $I = [e_1, \hdots, e_k]$ be an ordered subset of
$\mathcal{O}$, let $\alpha = [\alpha_1, \hdots, \alpha_k] \in
\{+,-\}^k$  and let $\mathcal{O}' = \mathcal{O}[I^{\alpha}] =
(\mathcal{O} \cup p)$ be the lexicographic extension of $\mathcal{O}$.

\begin{prop}\label{prop:IndexPUnequalZero}
The index $\ind_Y$ wrt.\ $I$ for a cocircuit $Y$ in $\mathcal{O}[I^{\alpha}]$ with $Y_p \neq 0$  is always at most $k$ and 
$Y_{e_{\ind_Y}} = Y_p$.  
\end{prop}

\begin{proof}
Because $Y_p \neq 0$ and Proposition \ref{prop:uniqueEdges2}, (ii) we obtain that $Y \setminus p$ is a cocircuit in $\mathcal{O}$.
The rest follows from the definition of the lexicographic extension.
\end{proof}

Recall that we call an extension {\em principal} if we add a point to a flat in general position in the matroid.  

\begin{prop}\label{prop:lexExtensionIsFree}
  The lexicographic extension $\mathcal{O}[I^{\alpha}]$ is a principal extension of $\mathcal{O}$ adding a point to the flat $\cl(e_1,
  \hdots, e_k)$.
\end{prop}

\begin{proof}
For each covector $X$ in $\mathcal{O}$ and its extension $X'$ in  $\mathcal{O}[I^{\alpha}]$ holds 
\[p \in z(X') \text{ iff } \{e_1, \hdots, e_k\} \subseteq z(X) \text{
  hence }\cl(e_1, \hdots, e_k) \subseteq z(X). \]
Assume $p$ were not in general position in $\cl(e_1,
  \hdots, e_k)$. Then $p$ is contained in a non-spanning cycle in that flat. Hence, there exists a cocircuit avoiding that cycle, which is non-zero on $p$, a contradiction.
\end{proof}

\begin{prop}\label{prop:thirdSmallProp}
Let $Y$ be a cocircuit in $\mathcal{O}[I^{\alpha}]$ with $Y_p = 0$. 
If the index $\ind_Y$ wrt. $I$ of $Y$ is $k+1$ then $Y$ is an old cocircuit otherwise it is a new cocircuit.
\end{prop}

\begin{proof}
  It is clear that if $\ind_Y < k+1$ and $Y$ were an old
  cocircuit, then $Y_p = 0$ would not hold.  If $\ind_Y = k+1$ and $Y$
  were a new cocircuit, it would have been derived from an edge $X
  \circ Z$. But then $I \subseteq z(X \circ Z)$ would imply $\rk(z(Y)) = \rk(z(Y \setminus p))$ because $p \in \cl(I)$, and hence $Y \setminus p$ would be a cocircuit already in $\mathcal{O}$. 
\end{proof}

We come to a technical lemma.

\begin{lemma}\label{lem:corrCocircuit}
  Let $Y$ be a new cocircuit in $\mathcal{O}[I^{\alpha}]$. Then we necessarily have $Y = {\left((X \setminus p) \circ (Z \setminus p), 0\right)}$ for two unique old cocircuits $X$ and $Z$ in $\mathcal{O}[I^{\alpha}]$
  such that $(X \setminus p) \circ (Z \setminus p)$ is an edge in
  $\mathcal{O}$ and $\sep(X,Z) = \{p\}$. The indices of $X$ and $Z$
  differ and are smaller than $k+1$. Let $i = \ind_X < j = \ind_Z$.
  We have
\[0 \neq X_{e_i} = Y_{e_i} = -Z_{e_j} = -Y_{e_j} \text{ and } \ind_Y = \ind_X < k.  \] 
\end{lemma}

\begin{proof}
From Proposition \ref{prop:thirdSmallProp} follows that $\ind_Y < k+1$. The first statement is Theorem \ref{theorem:ext}. 
Because $X \setminus p$ and $Y \setminus p$ are conformal and $X_p = -Z_p$, they must have different indices. 
The rest is immediate.
\end{proof} 

\begin{remark}\label{rem:corrCocircuit} We call the
    cocircuit $Z$ in Lemma \ref{lem:corrCocircuit} with higher index the cocircuit {\em corresponding to} $Y$.  \end{remark} 

  We can always extend an oriented matroid program to a bounded one.

  \begin{lemma}\label{lem:LexExtensionIsBounded} Let
    $(\mathcal{O},g,f)$ be an oriented matroid program.  Let
    $\mathcal{O}' = \mathcal{O} \cup p = \mathcal{O}[I]^+$ be a
    positive lexicographic extension of $\mathcal{O}$ in general
    position where $I$ is a base of $\mathcal{O} \setminus f$. Then
    $(\mathcal{O}',p,f)$ is bounded. It is feasible if and only if
    $(\mathcal{O},g,f)$ is.  \end{lemma} 

\begin{proof} Let
    $\mathcal{O}$ be of rank $\rk$. Then $I$ $= [e_1, \hdots, e_r]$ is
    an ordered set of independent elements in $\mathcal{O} \setminus
    f$.  If $(\mathcal{O},g,f)$ is feasible we have a cocircuit $X
    \neq 0$ with $X_e \in \{0,+\}$ for all $e \in E$ and also $X_p \in
    \{0,+\}$ since the lexicographic extension is
      positive. The program remains feasible.  Because $I$ is a base
    in $\mathcal{O}$ Proposition \ref{prop:thirdSmallProp} yields that
    any cocircuit $X$ with $X_p = 0$ must be new.  Lemma
    \ref{lem:corrCocircuit} yields two indices $i \neq j \leq \rk$ with
    $X_{e_i} = - X_{e_j} \neq 0$ hence $X_e = -$ for an $e \in E$ and
    $X$ is not in the feasible region. The program is bounded.
  \end{proof} 
\section{Proof of Theorem~\ref{cor:lexExtStaysEucl2a}}\label{sec:three}
  \begin{thmn}[\ref{cor:lexExtStaysEucl2a}] Let $\mathcal{O}$ be a
    Euclidean oriented matroid and let $\mathcal{O} \cup p$ be a
    lexicographic extension of $\mathcal{O}$.  Then $(\mathcal{O} \cup
    p,p,f)$ is a Euclidean matroid program.  \end{thmn} Before we go 
    into the details, we give a brief sketch of the proof.
  \subsection{Sketch of proof} We assume that there is a directed
  cycle $\mathcal{C}$ in $(\mathcal{O} \cup p,p,f)$ and show that
  there is already a directed cycle in an oriented matroid program in
  $\mathcal{O}$.  First, all cocircuits in $\mathcal{C}$ are old
  because $p=+$ holds for them, see Proposition
  \ref{prop:uniqueEdges2}, and have the same f-value, see Proposition
  \ref{prop:PathGraphWellKnownFacts}.  We show in Subsection
  \ref{subs:dirOfEdgesOfOld} that the edges between these old
  cocircuits are already edges in $\mathcal{O}$ and in Subsection
  \ref{section:DirGraphpf} that the directions of these edges are
  always preserved in some oriented matroid programs of $\mathcal{O}$.
  Hence, we always obtain a directed cycle in a matroid program of
  $\mathcal{O}$.  
  
  \subsection{Directions of edges of old cocircuits in
    single-element extensions}\label{subs:dirOfEdgesOfOld} We begin by
  examining the directions of old edges.  
    \begin{prop}\label{prop:EdgeStaysEdge} Let $(\mathcal{O},g,f)$ be
      an oriented matroid program and let $\mathcal{O}' = \mathcal{O}
      \cup p$ be a non-trivial single-element extension of
      $\mathcal{O}$.  Let $X,Y$ be two old conformal, comodular
      cocircuits in $(\mathcal{O}',g)$.  Then $X \setminus p,Y
      \setminus p$ are conformal cocircuits in $(\mathcal{O},g)$.
      They are also comodular in $(\mathcal{O},g)$ if $(X \circ Y)_p
      \ne 0$ or if $\mathcal{O} \cup p$ is a principal extension (or
      if both cases hold).  \end{prop} 

\begin{proof} Clearly, $X
    \setminus p, Y \setminus p$ are conformal cocircuits in
    $(\mathcal{O},g)$.  If $(X \circ Y)_p \neq 0$ then $z(X \circ Y) =
    z(X \circ Y) \setminus p$ hence $((X \circ Y) \setminus p)$ is an
    edge in $\mathcal{O}$.  If $\mathcal{O} \cup p$ is a principal
    extension, we have $p \in \cl_{\mathcal{M}'}(F \setminus p)$ for a
    flat $F \setminus p$ and all its supersets in $\mathcal{M}$.
    Hence if $p \in \cl_{\mathcal{M}'}(z(X \setminus p))$ and $p \in
    \cl_{\mathcal{M}'}(z(Y \setminus p))$ we have $F \setminus p
    \subseteq z(X \setminus p)$ and $F \setminus p \subseteq z(Y
    \setminus p)$ in $\mathcal{M}$.  We obtain $F \setminus p
    \subseteq z(X \setminus p) \cap z(Y \setminus p) = z((X \circ Y) \setminus p)$, therefore $p \in
    \cl_{\mathcal{M}'}(z((X \circ Y)) \setminus p)$ and $((X \circ Y)
    \setminus p)$ is an edge in $\mathcal{O}$.  \end{proof}

  \begin{remark}  If $(X \circ Y)_p = 0$ and
      $\mathcal{O} \cup p$ is not a principal extension it is possible
      that $X \setminus p$ and $Y \setminus p$ are not comodular in
      $\mathcal{O}$ e.g. if $\mathcal{O}$ has rank 3 and $z(X
      \setminus p)$ and $z(Y \setminus p)$ are disjoint lines in the
      underlying matroid intersecting at the point $p$ in the
      extension.  Because the lexicographic extension is principal
      (see Proposition \ref{prop:lexExtensionIsFree}), the two
      cocircuits are always comodular in $\mathcal{O}$ in that case.
     \end{remark} 

We introduce a new function to abbreviate our notions.  

\begin{definition}[Dir-function] Let
    $X,Y$ be cocircuits in \linebreak[3] $(\mathcal{O},g,f)$ and a comodular pair
    with $X_g = Y_g \neq 0$. Then we define
    \[\Dir_{(\mathcal{O},g,f)}(X,Y) = \left\{ \begin{array}[h]{lcr}
        + & \text{ iff \quad} X \rightarrow Y,\\
        - & \text{ iff \quad} X \leftarrow Y,\\
        0 & \text{ iff \quad} X \leftrightarrow Y.  \end{array}
    \right\} \text{ in } (\mathcal{O},g,f).\] If $X,Y$ are comodular
    cocircuits with $X_g \neq 0$ and $Y_g = 0$ (or vice versa) then we
    define $\Dir(X,Y) = Y_f$ (or $\Dir(X,Y) = -X_f$ resp.).
  \end{definition} 
 
Finally, we show that the directions of old edges are preserved in extensions.  

\begin{lemma}[direction preservation of
    old edges]\label{lem:DirPresOfOutsideEdges} Let
      $(\mathcal{O},g,f)$ be an oriented matroid program, and let
      $\mathcal{O}' = \mathcal{O} \cup p$ be a non-trivial
      single-element extension of $\mathcal{O}$.  Let $X,Y$ be two
    old cocircuits in $(\mathcal{O}',g)$ with $X \circ Y$ being an
    edge in $\mathcal{O}'$.  If $((X \circ Y) \setminus p)$ is also an
    edge in $\mathcal{O}$ then cocircuit elimination of $g$ between
    $-X$ and $Y$ yields an old cocircuit and it holds \[
    \Dir_{(\mathcal{O}',g,f)}(X,Y) =
      \Dir_{(\mathcal{O},g,f)}(X \setminus p, Y \setminus p). \]
  \end{lemma} 

\begin{proof} Since $F = z((X \setminus p) \circ (Y
    \setminus p))$ is a coline in $\mathcal{M}$, $\cl_{\mathcal{M}}(F
    \cup g)$ is a hyperplane in $\mathcal{M}$.  Then cocircuit
    elimination of $g$ between $-X \setminus p$ and $Y \setminus p$ in
    $\mathcal{O}$ yields a unique cocircuit $Z$ that can be extended
    to an old cocircuit $Z'$ in $\mathcal{O}'$. We have $z(Z') =
    cl_{\mathcal{M}'}(F \cup g)$.  Cocircuit elimination of $g$
    between $-X$ and $Y$ in $\mathcal{O}'$ yields a unique cocircuit
    $Z^1$ with $z(Z^1) \supseteq \cl_{\mathcal{M}'}(z(X \circ Y) \cup
    g) \supseteq z(Z')$.  For rank reasons, equality holds,
    so $Z^1 = \pm Z'$. Let $e \neq p$ with $e \in z(X \setminus p)$
    and $e \notin z(Y \setminus p)$ then we obtain $Z'_e = Z^1_e$
    hence $Z' = Z^1$ which proves the lemma.  \end{proof}
  \subsection{The directions of the cocircuit graph of the
    lexicographic extension $(\mathcal{O} \cup
    p,p,f)$}\label{section:DirGraphpf} We can describe completely the
  directions of the graph $G_f$ of $(\mathcal{O} \cup p,p,f)$ if we
  know the directions of some edges in the graphs $G_f$ of
  $\mathcal{O}$.  Until the end of this chapter, let $\mathcal{O}$ be
  an oriented matroid of rank $\rk$ with groundset $E$ and let $f \in
  E$.  Let ${I} = [e_1, \hdots, e_k]$ with $k \leq \rk$ be an ordered
  set of independent elements of $E$, where possibly $f \in I$.  Let
  $\mathcal{O}' = \mathcal{O} \cup p$ be the lexicographic extension
  $\mathcal{O}[e^+_1, \hdots, e^+_k]$.  

For the next
  technical lemmas, we fix some assumptions.

  \begin{definition}[Assumptions A] Let $X,Y$ be two conformal
    cocircuits with $X_p = Y_p \neq 0$ and $X \circ Y$ being an edge
    in $\mathcal{O}'$.  Let $P$ be the cocircuit yielded by
    elimination of $p$ between $-X$ and $Y$.  Let $\ind_X, \ind_Y$ be
    the indices of $X, Y$ wrt. $I$.  We assume $\ind_X \leq \ind_Y$.
    If $i = \ind_X = \ind_Y$ and thus $X_{e_i}=Y_{e_i}\ne 0$, let $ Z$
    denote the cocircuit yielded by elimination of $e_i$ between $-X$
    and $Y$.  If $\ind_X < \ind_Y$ let $Z = Y$.
  \end{definition} 

\begin{prop}\label{prop:cocircuitCompatible} Let
    assumptions A hold.  Then $P$ is old if and only if $P = Z$.  If
    $P$ is new, then $Z$ and $P$ are compatible.  \end{prop}

  \begin{proof} We have $z(Z)
    \supseteq z(X \circ Y) \cup e_i$ and $z(P) \supseteq z(X \circ Y)
    \cup p$.  Regarding the indices we obtain $\ind_Z > \ind_X$ and
    $\ind_P \geq \ind_X$.  Let $e \in z(X) \cap \supp(Y)$.  Then we
    must have $Y_e = P_e = Z_e \neq 0$ in all cases.  

    If $P_{e_i} = 0$
    Proposition \ref{prop:unqieCocircuitOnEndge} yields $P = \pm Z$
    hence $P = Z$.  Also $Z_p = 0$ implies $P = Z$ and $P$ is old. On
    the other hand, if $P$ is old, we must have $P_{e_i} = 0$ hence $P
    = Z$.  As by assumption $Y_p \neq 0$ we have that $\ind_X <
    \ind_Y$ implies $P \neq Y = Z$ hence $P$ must be new.

Now let $P$ be new and thus derived from two compatible comodular cocircuits $P^1, P^2$. 
We have $z(P^1 \circ P^2) = z(X \circ Y)$ thus by  Proposition~\ref{prop:uniqueEdges2} (i) and Lemma~\ref{lem:corrCocircuit}
we may assume that  
$\ind_X \le \ind_{P^2}<\ind_{P^1}.$ 
We conclude $P^1_{e_i} = 0$ and Proposition~\ref{prop:unqieCocircuitOnEndge} again yields $P^1 = \pm Z$.
But because $P^1 \neq \pm X$ we have $P^1_e \neq 0$ and thus $P^1_e=P_e$. Hence $P_e=Z_e$ implies $P^1=Z$, and hence $P$ and $Z$ are compatible. 
\end{proof}

\begin{lemma}\label{lem:case3lexext}
  Let assumptions A hold.  Let $1 \leq i = \ind_X = \ind_Y \leq k$ and
  $f \neq e_i$ and thus $X_{e_i} = Y_{e_i}\neq 0$.  If $X_f \neq Y_f$
  or $X_f = Y_f = 0$ or  $\Dir_{e_i,f}(X,Y) \neq 0$ then 
\[\Dir_{e_i,f}(X,Y) = \Dir_{p,f}(X,Y) \text{ and } \Dir_{f,e_i}(X,Y) = \Dir_{f,p}(X,Y). \]
\end{lemma}

\begin{proof}

If $X_f \neq Y_f$ or $X_f = Y_f = 0$
Proposition \ref{prop:basicDirectionProps} (iii) yields $\Dir_{e_i,f}(X,Y) = \Dir_{p,f}(X,Y)$.
If $\Dir_{e_i,f}(X,Y) \neq 0$, we may assume that $X_f = Y_f \neq 0$. Hence if $f \in  I$ then $f = e_j$ with $j > i$.
In any case Proposition \ref{prop:cocircuitCompatible} yields either $Z = P$ or $P$ is compatible to $Z$.
We have $Z_f \neq 0$, and therefore $P_f = Z_f$ implying the assertion. The second equation (if defined) follows from Lemma \ref{prop:ChangeOfFAndG1} or Corollary \ref{cor:ChangeOfFAndG2} and $X_{e_i} = Y_{e_i} = X_p = Y_p$.
\end{proof}

\begin{lemma}\label{lem:case4lexExt}
Let assumptions A hold. 
Let $1 \leq i = \ind_X = \ind_Y \leq k$. Hence $X_{e_i} = Y_{e_i} = X_p = Y_p \neq 0$. Let $X_f = Y_f  \neq 0$ and $\Dir_{e_i,f}(X,Y) = 0$ or $f = e_i$.
Then 
\begin{itemize}
\item[(i)] $\Dir_{e_i,e_j}(X,Y) = 0$ for all $j$ with $i  < j < k+1 \Rightarrow \Dir_{p,f}(X,Y) = 0$.  
\item[(ii)] Otherwise $\Dir_{f,p}(X,Y) = \Dir_{e_i,e_j}(X,Y)$, where $j > i$ is the first index with
$\Dir_{e_i,e_j}(X,Y) \neq 0$.
\item[(iii)] It holds $\Dir_{f,e_j}(X,Y) = \Dir_{e_i,e_j}(X,Y)$.
\end{itemize}
\end{lemma}

\begin{proof}
$\Dir_{e_i,f}(X,Y) = 0$ and $f = e_i$ both imply $Z_f =  0$ and $\ind_Z = j > i$. 
Because $X_f = Y_f \neq 0$ we have $f \notin \{e_1, \hdots, e_{i-1}\}$.  

ad (i): 
We obtain $Z_{e_j} = 0$ for all $i < j < k+1$ hence
Proposition \ref{prop:IndexPUnequalZero} yields $Z_p = 0$ and Proposition \ref{prop:cocircuitCompatible} yields $Z = P$ hence $Z_f = P_f = 0$. 

ad (ii): We have $Z_{e_j} \neq 0$.
Lemma \ref{lem:DirPresOfOutsideEdges} yields that $Z$ is an old cocircuit. We obtain $Z_p = Z_{e_j}$  because $j > i$ is the first index where $\Dir_{e_i,e_j}(X,Y) \neq 0$.
Cocircuit elimination of $f$ between $-X$ and $Y$ yields also the cocircuit $Z$ which proves (ii).

ad (iii): This is Proposition \ref{prop:TransitivityOfDir}.
%
\end{proof}

\begin{lemma}\label{lem:case6lexExt}
Let assumptions A hold and let $Y_f \neq 0$. 
Let $1 \leq i = \ind_X < \ind_Y = j \leq k$. Then 
\[\Dir_{p,f}(X,Y) = \Dir_{e_i,f}(X,Y) = Y_f. \]
\end{lemma}

\begin{proof}
Proposition \ref{prop:cocircuitCompatible} yields that $Y$ and $P$ are compatible hence $P_f = Y_f$.
\end{proof}

We denote some small facts concerning cycles in $(\mathcal{O}',p,f)$ 
where $\mathcal{O}' = \mathcal{O} \cup p$ is a non-trivial single-element extension of $\mathcal{O}$.

 \begin{prop}\label{prop:cocircuitsStayCocircuitsInOriginal}
All cocircuits in a cycle in $(\mathcal{O}',p,f)$ are also cocircuits in $\mathcal{O}$
and all edges remain edges.
\end{prop}

\begin{proof}
This follows immediately from Proposition \ref{prop:uniqueEdges2} and Proposition \ref{prop:EdgeStaysEdge}. 
\end{proof}

For the last three lemmas, let $\mathcal{O}' = \mathcal{O} \cup p = 
 \mathcal{O}[e^+_1, \hdots, e^+_k]$ be a positive lexicographic extension of $\mathcal{O}$. 

 \begin{lemma}\label{lem:directedCycleNoIndex1}
All cocircuits in a directed cycle in $(\mathcal{O}',p,f)$ have the same index. 
\end{lemma}

\begin{proof}
  Because of Proposition \ref{prop:PathGraphWellKnownFacts} (iv) all
  cocircuits in the cycle must have the same f-value $\neq 0$.  Suppose that there is a directed cycle $X_0,\ldots X_n=X_0$ in which not all vertices have the same index. Then there exist $X_i, X_j$ such that $\ind_{X_i} < \ind_{X_{i+1}}$ and $\ind_{X_j} > \ind_{X_{j+1}}$.  But then Lemma \ref{lem:case6lexExt} implies $\Dir_{p,f}(X_i,X_{i+1}) = -\Dir_{p,f}(X_j,X_{j+1})$ contradicting the cycle being directed.
\end{proof}

 \begin{lemma}\label{lem:case6lexExt2}
Let $\mathcal{G}$ be a directed cycle in $(\mathcal{O}',p,f)$ where all cocircuits of $\mathcal{G}$ have the same index $i$.
If one edge of the cycle is directed in $(\mathcal{O}',e_i,f)$, we have a directed cycle in $(\mathcal{O}',e_i,f)$.
\end{lemma}

\begin{proof}
Lemma \ref{lem:case3lexext} yields that the edge is directed as in $(\mathcal{O}',p,f)$. All other edges are either undirected or
directed as in $(\mathcal{O}',p,f)$. We obtain a directed cycle in $(\mathcal{O}',e_i,f)$.
\end{proof}

 \begin{lemma}\label{lem:DirectedCycleInEj}
Let $\mathcal{G}$ be a directed cycle in $(\mathcal{O}',p,f)$ where all cocircuits of $\mathcal{G}$ have the same index $i$.
If all edges of the cycle are undirected in $(\mathcal{O}',e_i,f)$ we have a directed cycle in $(\mathcal{O}', e_i,e_j)$ (and in $(\mathcal{O}', f,e_j)$) for a $j > i$.
This case cannot occur if $f$ is in general position and if $f \neq e_i$.
\end{lemma}

\begin{proof}
Proposition \ref{prop:DirectedCycleYieldsDirectedCycle} yields a directed cycle $\mathcal{G}' = \pm \mathcal{G}$ in $(\mathcal{O}',f,p)$.
  If all edges of $\mathcal{G}'$ were undirected in
  $(\mathcal{O}',e_i,e_j)$, Lemma \ref{lem:case4lexExt} would yield
  that the edges are also undirected in $(\mathcal{O}',f,p)$
  contradicting the assumptions.  Hence, we obtain at least one
  directed edge in $(\mathcal{O}',e_i,e_j)$ for some $j > i$. We take
  the smallest index $j > i$ where such an edge exists.  Lemma
  \ref{lem:case4lexExt} yields that the edge is directed as in
  $(\mathcal{O}',f,p)$. If we have another directed edge in
  $(\mathcal{O}',e_i,e_j)$ Lemma \ref{lem:case4lexExt} yields again
  (because the edge is undirected in all $(\mathcal{O}',e_i,e_k)$ with
  $k < j$) that it must be directed like in $(\mathcal{O}',f,p)$.  That means that in
  $(\mathcal{O}',e_i,e_j)$ all edges of the cycle are undirected or
  directed as in $(\mathcal{O}',f,p) $.
  We obtain a directed cycle in $(\mathcal{O}',e_i,e_j)$. 
  Because of Lemma \ref{lem:case4lexExt}  (iii), the proof also works if we substitute $e_i$ by $f$ everywhere.
  The second
  statement is immediate.
\end{proof}

\subsection{First Main Theorems}

We formulate the main theorem with the least possible assumptions.
$\mathcal{O}$ does not need to be fully Euclidean, only a few programs
of $\mathcal{O}$ and here even only the contractions of these
programs.  (Recall that a contraction of a Euclidean oriented matroid
program remains Euclidean, see \cite{2},10.5.6.)

\begin{theorem}\label{theo:lexExtStaysEucl2}
Let $\mathcal{O}$ be an oriented matroid of rank $\rk$ with groundset $E$ and let $f \in E$. 
Let ${I} = [e_1, \hdots, e_k]$ with $k \leq \rk$ be an ordered set of independent elements of $E$.
Let $l  \leq  k$ be the smallest index
such that $f \in \cl(\{e_1, \hdots, e_l\})$. If there is no such $l$ let $l = k+1$ (in that case $k+1 \le \rk$).

\begin{enumerate}
\item Let  $(\mathcal{O},e_1,f), (\mathcal{O} / \{e_1\},e_2,f), \hdots, (\mathcal{O} / \{e_1, \hdots, e_{l-2}\},e_{l-1},f)$ be Euclidean oriented matroid programs and
\item If $f$ is not in general position or if $f \in {I}$ let  $(\mathcal{O} / \{e_1, \hdots, e_{i-1}\},e_i,e_j)$ be Euclidean oriented matroid programs for all $1 \le i \le min(l,k-1)$ and $i < j \le k$.
\end{enumerate}
Let $\mathcal{O}' = \mathcal{O} \cup p$ be the lexicographic extension $\mathcal{O}[e^+_1, \hdots, e^+_k]$.  
Then $(\mathcal{O}',p,f)$ is a Euclidean oriented matroid program.
\end{theorem}

\begin{proof}
We show that the graph $G_f$ of $(\mathcal{O}',p,f)$ has no directed cycles.
Proposition \ref{prop:PathGraphWellKnownFacts} (iv) shows that the
cocircuits of a directed cycle must all have the same $f$-value $v_f
\neq 0$.  Lemma \ref{lem:directedCycleNoIndex1} shows that the
cocircuits of the cycle all have the same index $i$. We obtain that $f
\notin \cl\{e_1, \hdots,e_{i-1} \}$ hence $i \le l$.  If one of the edges of the cycle is directed in $(\mathcal{O}',e_i,f)$, Lemma \ref{lem:case6lexExt2} yields
a directed cycle in $(\mathcal{O}',e_i,f)$.  Otherwise (and that case
can only happen if $f$ is not in general position or if $f = e_i$),
Lemma \ref{lem:DirectedCycleInEj} yields a directed cycle in
$(\mathcal{O}', e_i,e_j)$ {for some $j > i$}.  Lemma
\ref{lem:DirPresOfOutsideEdges} yields that we obtain directed cycles
in either $(\mathcal{O},e_i,f)$ or $(\mathcal{O},e_i,e_j)$.  Since all
indices of the cocircuits {of these cycles} are $i$, in the
contractions $(\mathcal{O}/\{e_1, \hdots, e_{i-1}\},e_i,f)$ and
$(\mathcal{O}/\{e_1, \hdots, e_{i-1}\}, e_i, e_j)$ the cocircuits stay
cocircuits and we obtain a directed cycle in at least one of them,
contradicting the assumptions.
\end{proof}
\begin{remark}\label{rem:SubstitutionMainTheoremAssumptions}
In assumption (ii) of the Theorem we can substitute $e_i$ by $f$, the proof works analogously. 
\end{remark}

We get immediately:

\begin{corollary}\label{cor:lexExtStaysEucl2}
Let $\mathcal{O}$ be an oriented matroid of rank $\rk$ with groundset $E$ and let $f \in E$ be in general position. 
Let ${I} = [e_1, \hdots, e_{k}]$ with $k \le \rk$ be an ordered set of independent elements of $E \setminus f$ such that
\[(\mathcal{O},e_1,f), (\mathcal{O} / \{e_1\},e_2,f), \hdots, (\mathcal{O} / \{e_1, \hdots, e_{k-1}\},e_k,f)\] are Euclidean oriented matroid programs.
Let $\mathcal{O}' = \mathcal{O} \cup p$ be the lexicographic extension $\mathcal{O}[e^+_1, \hdots, e^+_k]$.  
Then $(\mathcal{O}',p,f)$ is a Euclidean oriented matroid program.
\end{corollary}

 \section{Proof of Theorem \ref{theorem:SecondMainTheorem}}\label{sec:four}
 
\begin{thmn}[\ref{theorem:SecondMainTheorem}]
A lexicographic extension of a Euclidean oriented matroid
is Euclidean.
\end{thmn}
 
 \subsection{Sketch of Proof}

 To show that the program $(\mathcal{O}' = \mathcal{O} \cup p,g,f)$ of
 a lexicographic extended oriented matroid remains Euclidean if
 $\mathcal{O}$ is Euclidean is slightly more complicated than the proof of our
 first theorem, where $p$ was a new hyperplane at infinity.  {Like before we show that a directed
   cycle $P$ in the program $(\mathcal{O}',g,f)$ yields a directed
   cycle in a program of $\mathcal{O}$.  We know already from Lemma
   \ref{lem:DirPresOfOutsideEdges} that the directions of edges
   between old cocircuits in $(\mathcal{O}',g,f)$ are preserved in
   $(\mathcal{O},g,f)$ hence we may assume that $P$ contains new
   edges.  These new edges do not exist in $\mathcal{O}$ we have to
   substitute them with old edges preserving their directions.  We
   apply two approaches to achieve that.  On the one hand we prove a
   Lemma in Subsection \ref{subs:projTriangleLemma}, the {\em
     Projection Lemma} \ref{lem:ProjectionLemma}.  so that we can
   project new edges (and their directions) to edges in the program
   $(\mathcal{O} \cup p,p,f)$ obtaining a projected directed cycle
   there.  On the other hand we examine all constellations of new
   edges together with their corresponding cocircuits.  In some
   constellations corresponding edges appear and we obtain a
   corresponding cycle.  We examine the directions of these edges in
   Subsection \ref{section:DirGraphgf} again frequently using the Projection
   Lemma.  The directions of the edges that give rise to  the new
   cocircuits of $P$ (analyzed in Subsection
   \ref{subsection:dirOfEdgesWithNewCocInSEEctensions}) do matter as well. It
   is important that they do not suddenly  change since  otherwise we lose
   information about the direction of the corresponding edges.  We
   call cycles / paths with that property {\em normalized} and prove in
   Subsection \ref{subs:normalizedPath} that we can always get a
   normalized version of a directed cycle.  We use the {\em Triangle
     Lemma} \ref{lem:TriangleLemma} and the {\em Zero Line Lemma}
   \ref{lem:ZeroLineLemma} proven in Subsection
   \ref{subs:projTriangleLemma} here.  Hence we may assume that our
   directed cycle $P$ is normalized.  In Subsection
   \ref{section:DirGraphgf} we show that if $P$ has old cocircuits or
   cocircuits with different indices we obtain a corresponding cycle
   $Q$ that is directed in $(\mathcal{O},g,f) $.  Finally assuming $P$
   having only new cocircuits with index $i$ the corresponding cycle
   (if existing) is either directed like $P$ in $(\mathcal{O},g,f)$ or
   in $(\mathcal{O},g,e_j)$ for a $j > i$, or is undirected. But if the
   corresponing cycle is undirected or even not existing we use our
   first approach and project $P$ to a directed cycle in $(\mathcal{O}
   \cup p, p,f)$ and are done.  We show the final result in Subsection
   \ref{subsection:secondMainResult} and a Corollary with less
   assumptions for the case that $g$ and $f$ are in general position.}
 
\subsection{The Projection Lemma and the Triangle Lemma}\label{subs:projTriangleLemma}

{The lemmas in this preparatory section are} results concerning edges in oriented matroids whose zero sets intersect in a coplane 
(a flat of corank $3$) in the underlying matroid. 
We call $(F,G,H)$ a {\em modular triple} of covectors in an oriented matroid 
if $(F,G)$, $(G,H)$ and $(F,H)$ are modular pairs in the underlying matroid.

The following is immediate.

\begin{prop}\label{prop:crossingPoint}
  Let $(F,G)$ be a pair of edges in an affine oriented matroid
  $(\mathcal{O},g)$ with $z(F \cap G)$ being a coplane in the
  underlying matroid $\mathcal{M}(\mathcal{O})$.  Then there is up to sign reversal a unique
  cocircuit $C$ lying on the line through  $F$ and on the line through $G$.  
\end{prop}

\begin{definition}\label{def:crossingPoint}
  Let the assumptions of  Proposition~\ref{prop:crossingPoint}  hold.
  If $C_g \neq 0$ we call the cocircuit $C$ with
  $C_g = +$ the {\em crossing point of $F$ and $G$}. If $C_g = 0$, the crossing point is arbitrarily chosen between $C$ and $-C$.
\end{definition}

We use the assumptions of the next proposition in the whole Section.

\begin{prop}\label{prop:additionalAssumption}
Let $(C,X,Y)$ be a modular triple of cocircuits in an oriented matroid of rank $\geq 3$ with $z(C \circ X \circ Y)$ being a coplane in the underlying matroid.
Then $(C \circ X,C \circ Y, X \circ Y)$ is a modular triple of edges with different zero sets, which means that $C,X,Y$ do not lie on one line. 
\end{prop}

\begin{proof}
Because $(C,X,Y)$ is a modular triple, $C \circ X, C \circ Y, X \circ Y$ are edges. 
Because $z(C \circ X \circ Y)$ is a coplane, it is clear that they form a modular triple, too.
All these edges have different zero sets; otherwise
if we had e.g. $z(C \circ X) = z(C \circ Y)$, we would obtain $z(C \circ X) = z(C \circ X \circ Y)$ and $z(C \circ X)$ would not be a coline. 
\end{proof}

\begin{prop}\label{prop:edgesBetweenCrossingEdges}
  Let $F$ and $G$ be two edges of an affine oriented matroid
  $(\mathcal{O},g)$ with $z(F \cap G)$ being a coplane in $\mathcal{M}(\mathcal{O})$.  Let $C$ be the crossing point of $F$ and $G$.
  Let $X$ resp. $Y$ be cocircuits in $(\mathcal{O},g)$ lying on the line through $F$ resp. through $G$ and $C \neq X
  \neq Y \neq C$. If $z(X) \cap z(Y) \cap \supp(C) \ne
  \emptyset$, then $X \circ Y$ is an edge, 
  $(C, X, Y)$ a modular triple with $z(C \circ X \circ Y)$ being a coplane in $\mathcal{M}(\mathcal{O})$.
\end{prop}

\begin{proof}
We have $z(F) = z(X \circ C) = z(C \circ X)$ and $z(G) = z(Y \circ C) = z(C \circ Y)$. With $H=F\circ G$ we obtain $z(F \cap G) = z(H) = z(C \circ X \circ Y)$ being a coplane.
Since $z(H) \subsetneq z(X \circ Y)$ we obtain $X \circ Y$ being an edge, and hence $(C,X,Y)$ being a modular triple.
\end{proof}

The proofs of the following lemmata are all reduced to the case of pseudoline arrangements (oriented matroids of rank $3$).
First, recall that oriented matroid programs of rank $3$ are always Euclidean, see  \cite{2}, Proposition 10.5.7.
The next lemma, called the 'Projection Lemma' yields a connection between oriented matroid programs with exchanged target functions. 

\begin{figure}
  \centering
\includegraphics[scale=0.7]{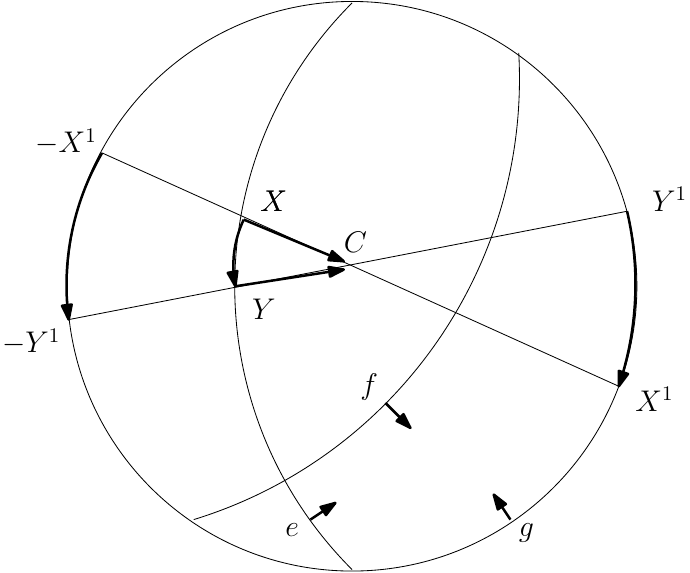} 
  \caption{Projection Lemma \label{fig:projectionLemma}}
\end{figure}

 \begin{lemma}[Projection Lemma]\label{lem:ProjectionLemma}
 Let $(C,X,Y)$ be a modular triple of cocircuits in an oriented matroid program $(\mathcal{O},g,f)$  of rank $\geq 3$ with $z(C \circ X \circ Y)$ being a coplane in $\mathcal{M}(\mathcal{O})$.
Let $C_g = X_g = Y_g = +$.
Let $e \in \supp(C)$ and $(X \circ Y)_e = 0$.
Let cocircuit elimination of $g$ between $-X$  and $C$ yield $X^1$
and between $-Y$  and $C$ yield $Y^1$. 
Then $X^1 \circ Y^1$ is an edge, and

\[ X \rightarrow_{g,f} Y  \iff Y^1 \rightarrow_{e,f} X^1 \iff -X^1 \rightarrow_{e,f} -Y^1. \]

\end{lemma}

\begin{proof}
  Since $z(X), z(Y), z(C), z(X^1), z(Y^1) \supset z(X \circ Y \circ C)
  = S$ we can continue with the contraction $P' = (\mathcal{O} / S,
  g,f)$ which is a rank $3$ oriented matroid program.  Necessarily we must have $X^1
  \neq Y^1$ for otherwise  $z(C \circ X) = z(C \circ X^1) =
  z(C \circ Y^1) = z(C \circ Y)$ contradicting the assumptions.
  $(X^1 \circ Y^1)_g = 0$ then implies that $X \circ Y$ is an edge. Cocircuit
  elimination of $g$ between $-X$ and $Y$ yields $Z$ with $Z_e = Z_g =
  0$.  Cocircuit elimination of $e$ between $-{Y}^1$
  and ${X}^1$ yields $Z^1$ with $Z^1_e = Z^1_g = 0$.
  This implies $Z=\pm Z^1$. Let $e_1 \in z({Y \circ C
  }) \setminus z( { X \circ C })$, then $X_{e_1}\ne 0$
  but $C_{e_1}=0$ and $X_{e_1}^1=-X_{e_1}$. Now $Z_{e_1}=-X_{e_1}$ and
  $Z^1_{e_1}=X^1_{e_1}$ implies $Z^1=Z$ implying the first
  equivalence. The second equivalence follows from Proposition
  \ref{prop:basicDirectionProps}.
\end{proof}

We formulate the so-called 'Triangle Lemma'.

 \begin{lemma}[Triangle Lemma]\label{lem:TriangleLemma}
  Let $(C,X,Y)$ be a modular triple of cocircuits in an oriented matroid program $(\mathcal{O},g,f)$  of rank $\geq 3$ with $z(C \circ X \circ Y)$ being a coplane in $\mathcal{M}(\mathcal{O})$.
Let $X_g = Y_g = +$ and $C_g = +$ or $C_g = 0 \neq C_f$.
Then the triangle $(X, Y, C)$ connected via the lines through $(X \circ Y, X \circ C, Y \circ C)$ is not a directed cycle. 
 \end{lemma}
 
 \begin{proof}
 We omit the obvious case $C_g = 0 \neq C_f$, so let $C_g = +$. 
First, we assume $f \in z(C \circ X \circ Y) = S$. Then we have $X_f = Y_f = C_f = 0$ and (using Proposition \ref{prop:PathGraphWellKnownFacts}(ii))
 all the edges of the triangle are undirected.
 Otherwise, and if the triangle was a directed cycle, we would obtain a directed cycle in the contraction $P' = (\mathcal{O} / S, g,f)$. 
 But that is impossible because $P'$ is of rank 3 and hence Euclidean.
 \end{proof}

If we have $C_g = C_f = 0$ for the crossing point $C$ of two edges, neither the triangle nor the projection lemma yield any information here. We show another lemma regarding that case.

 \begin{lemma}[Zero-Line Lemma]\label{lem:ZeroLineLemma}
     Let $(C,X,Y)$ be a modular triple of cocircuits in an oriented matroid program $(\mathcal{O},g,f)$  of rank $\geq 3$ with $z(C \circ X \circ Y)$ being a coplane in $\mathcal{M}(\mathcal{O})$. Let $X_g = Y_g = +$ and let $C_g = C_f = 0$.  Let
   $X^1, Y^1$ be different cocircuits in $(\mathcal{O},g)$ lying on the line through $F =  X \circ C$ resp. $G = Y \circ C$  with $X \neq X^1 \neq C$ and  $Y \neq Y^1 \neq C$ 
   with $X^1 \circ Y^1$ being an edge.
   Then $X^1
   \leftrightarrow X$ and $Y^1 \leftrightarrow Y$ and furthermore
\[ X^1 \rightarrow Y^1 \Leftrightarrow X \rightarrow Y. \]
 \end{lemma}
 
 \begin{proof} 
    $C_f = 0$ implies that all edges
   along the line through $F$ and along the line through $G$ are undirected.  
   If $X^1 \rightarrow Y^1$ and $X \rightarrow Y$ would not hold, then Proposition \ref{prop:PathGraphWellKnownFacts}(ii)
   would yield a directed cycle $(X^1, \hdots, Y^1, \hdots, Y, \hdots, X, \hdots, X^1)$  in $P$ and also
   in the contraction $P' = (\mathcal{O} / (z(F) \cap z(G)),g,f)$
   contradicting the Euclideaness of $P'$.
 \end{proof} 
 
 \subsection{Directions of Edges with New Cocircuits in Single-element Extensions}\label{subsection:dirOfEdgesWithNewCocInSEEctensions} 
{In this section let $(\mathcal{O},g,f)$ be an oriented matroid program of rank $\rk$ and 
let $\mathcal{O}' = \mathcal{O} \cup p$ be a non-trivial single-element extension of $\mathcal{O}$.
Even for the case of an arbitrary single element extension we can show that edges between old and new cocircuits behave in an obvious way.  }

\begin{prop}\label{prop:derivedDirectedLikeOrig}
Let $X$ be a new cocircuit \textcolor{black}{in $\mathcal{O}'$} with $X_g = +$ derived from the two cocircuits 
$X^1,X^2$ with $X^1_g = X^2_g = +$. Then we have
\[ \textcolor{black}{\Dir_{(\mathcal{O}',g,f)}(X^1,X) = \Dir_{(\mathcal{O}',g,f)}(X,X^2) =  \Dir_{(\mathcal{O}',g,f)}(X^1,X^2).} \]
\end{prop}

\begin{proof}
Proposition \ref{lem:DirPresOfOutsideEdges} yields \textcolor{black}{$\Dir_{(\mathcal{O},g,f)} (X^1,X^2) =  \Dir_{(\mathcal{O}',g,f)}(X^1, X^2)$}
for the old cocircuits $X^1,X^2$ and the rest follows from Proposition \ref{prop:basicDirectionProps} (iv).
\end{proof}

\begin{prop}\label{prop:derivedDirectedLikeOrigg0}
Let $X$ be an new cocircuit \textcolor{black}{in $\mathcal{O}'$} with $X_g = +$ derived from the two cocircuits 
$X^1,X^2$ with $X^1_g = +$ and $X^2_g = 0$. Then it holds 
\[ \Dir_{(\mathcal{O}',g,f)}(X^1,X) = X^2_f =  \Dir_{(\mathcal{O},g,f)}(X^1,X^2).\]
\end{prop}

\begin{proof}
Cocircuit elimination of $g$ between $-X^1$ and $X$ yields a cocircuit $W$ with $W_g = 0, W_p = X^2_p$.
It is clear that $W = X^2$.
\end{proof}

Now we examine the directions of new edges. Because we still work with arbitrary extensions,
we can not say very much in that case, but 
at least we can show two helpful propositions.

\begin{prop}\label{prop:CuttingConditionInExtension}
 Let $X,Y$ be conformal new cocircuits in
  $(\mathcal{O}',g)$ derived from the cocircuits $X^1, X^2$ and $Y^1, Y^2$
  with $X \circ Y$ being an edge. Let $F = X^1 \circ X^2$ and $G =  Y^1
  \circ Y^2$. Then there is a crossing point $C$ different from $X$ or $Y$ between the lines through $F$ and $G$  
  with $C_p \neq 0$ and $(C,X,Y)$ is a modular triple with $z(C \circ X \circ Y)$ being a coplane in $\mathcal{M}(\mathcal{O})$.
\end{prop}

\begin{proof}
We have $z(X) = z(F) \cup p$ and $z(Y) = z(G) \cup p$ hence $z(X \circ Y) = (z(F) \cap z(G)) \cup p $
\textcolor{black}{and $z(X) \neq z(Y)$ yields $z(F)\neq z(G)$.}
Because of $\rk(z(X \circ Y)) = \rk - 2$ we have $\rk(z(F) \cap z(G)) \geq r-3$ and equality follows from
$z(F) \neq z(G)$. Proposition \ref{prop:crossingPoint} yields the existence of a crossing point $C$.
We have $X \neq C \neq Y$
otherwise e.g.  $C = X$ would yield $z(X) = z(C) \supseteq z(G) \cup p = z(Y)$ hence $X = \pm Y$ which is impossible.
It follows $C_p \neq 0$ and Proposition \ref{prop:edgesBetweenCrossingEdges} yields the rest.
\end{proof}

That means that the Triangle Lemma, the Projection Lemma and the Zero-line Lemma can be applied in that case.
Applying the Projection Lemma \ref{lem:ProjectionLemma}
immediately yields a projection of undirected edges in another oriented matroid program. 

\begin{prop}\label{prop:projectionLemmaWithUndirectedEdges}
Let $X,Y$ be conformal new cocircuits in $(\mathcal{O}',g)$ derived from the cocircuits $X^1, X^2$ and $Y^1, Y^2$ with $X \circ Y$ being an edge. 
Let $C$ be the crossing point of $X^1 \circ X^2$ and $Y^1 \circ Y^2$ and let $C_g = +$.
Let cocircuit elimination of $g$ between $-C$ and $X$ yield $X^3$ and
between $-C$ and $Y$ yield $Y^3$. Then $X^3 \circ Y^3$ is an edge and it holds
\[ X \leftrightarrow _{g,f} Y\iff  X^3 \leftrightarrow_{p,f} Y^3 \]
\end{prop}

\subsection{Paths and Cycles in Single Element Extensions}\label{subs:normalizedPath}


\textcolor{black}{In this chapter we define the notion of a {\em normalized cycle} and prove that from each directed cycle we can derive a normalized version.}
  In the following if we have two comodular cocircuits $X,Y$ we say that {\em $X$ and $Y$ span a line} and with {\em the path from $X$ to $Y$}
  we always mean the shorter path from $X$ to $Y$ along the line through $X \circ Y$. 
 We start with a technical proposition. \textcolor{black}{$\mathcal{O}$ and $\mathcal{O}'$ are defined like in the previous subsection.} 

\begin{prop}\label{prop:QuadrangleConstellationNoCrossingPoint}
Let $X,Y$ be conformal new cocircuits in $\mathcal{O}'$ derived from the cocircuits $X^1, X^2$ and $Y^1, Y^2$ with $X \circ Y$ being an edge. Let $F = X^1 \circ X^2$ and $G = Y^1 \circ Y^2$ and let $F \neq G$. 
Let $C$ be the crossing point of $F$ and $G$. 
Let $C \neq X \neq Y \neq C$ and let 
$X^1_p =Y^1_p \neq 0$.
Then $C = \pm X^1$ if and only if $Y^1 = X^1$. 
\end{prop}

\begin{proof}
  It is clear that if $X^1 = Y^1$ implies $X^1 = Y^1 = \pm C$. For the
  other implication, it suffices to show that $C$
  and $Y$ are conformal. 
Assume there is $q\in
  \sep(C,Y)$. Since $X$ is conformal with $C$ as well as with $Y$, we
  must have $q \not \in \supp(X)$. Clearly, $q\ne p$, implying $q\in
  z(F) \subseteq z(C)$, contradicting $q \in \sep(C,Y)$.
\end{proof}

The next lemma says that in some constellations we can substitute an
edge between two new conformal cocircuits with a path of old
cocircuits having the same direction.
 
\begin{lemma}\label{lem:path1lem}
Let $X,Y$ be conformal cocircuits in $(\mathcal{O}',g)$ with $X_p = Y_p = 0$
and $X \circ Y$ being an edge.
Let $Y$ be new, derived from $Y^1 \circ Y^2$. 
If $X$ is new, let it be derived from $X^1 \circ X^2$ 
and let $X^1_p = Y^1_p$.
Then if either
\begin{enumerate}
\item $\Dir(Y^1,Y^2) = 0$ or
\item both cocircuits are new and $\Dir(X^1,X^2) = - \Dir(Y^1,Y^2) \neq 0$ 
\end{enumerate}
there is a path
\[ P = \begin{cases}
(X, \hdots, Y^i,Y) &\text{ if } X \text{ is old,}  \\
(X,X^i, \hdots, Y^i,Y) &\text{ if } X \text{ is new (with } i \in \{1,2\})  \\
\end{cases} \]
in $(\mathcal{O}',g)$ directed like the edge $X \circ Y$. All cocircuits of $P \setminus \{X,Y\}$ are old.
\end{lemma}

\begin{proof}
  If $X$ is new let $C$ be the crossing point between the edges $X^1 \circ X$ and $Y^1 \circ Y$.
  Proposition \ref{prop:CuttingConditionInExtension} yields $(C,X,Y)$ being a modular triple with $z(C \circ X \circ Y)$ being a coplane in $\mathcal{M}$.
  If $X$ is old let $e \in \supp(Y)\setminus \supp(X)$ and let $C$ be a cocircuit on the
  line through $Y^1 \circ Y$ with $C_e = 0$ and $C_g \neq -$.
   We have $z(C \circ Y \circ X) = z(Y^1 \circ Y \circ X) = z(Y \circ X) \setminus p$ hence $\rk(z(C \circ Y \circ X)) \geq \rk - 3$. 
    Equality holds because $p \notin z(C \circ Y)$ and $p \in z(X \circ Y)$. Hence $z(C \circ X \circ Y)$ is a coplane in $\mathcal{M}$.
    Applying Proposition \ref{prop:edgesBetweenCrossingEdges} with $Y$
  being the crossing point yields that $X$ and $C$ span a line and $(C,X,Y)$ is a modular triple.

  First we assume $C_g = +$.  Let $Q^1$ be the path on the line from
  $X$ to $C$, let $Q^2$ be the path on the line from $C$ to $Y$ and
  let $P = Q^1 \cup (Q^2 \setminus C)$.  Since $Y_g = X_g = C_g = +$
  and because $C_p \neq 0$ all cocircuits of $P$ except $X$ and $Y$
  are old and lie in $(\mathcal{O}',g)$.  In case (ii) $\Dir(Y_1,Y_2) =
  - Dir_(X_1,X_2) \neq 0$ implies $\Dir(X, C) = -\Dir(Y,C) = \Dir(C,Y)$.
  Hence $Q^1$ is directed like $Q^2$ and $P$ is a directed path.  It
  must be directed like the edge $X \circ Y$ otherwise we would obtain
  a directed cycle contradicting the Triangle Lemma
  \ref{lem:TriangleLemma}.  In case (i) the path $Q^2$ is undirected
  but then $Q^1$ must be directed like the edge $X \circ Y$ otherwise
  we obtain a directed cycle again.  In all cases $P$ is directed like
  the edge $X \circ Y$ and ends with $(\hdots, Y^i,Y)$ and starts with
  $(X, X^i, \hdots)$ if $X$ is new with $i \in \{ 1,2\}$ hence is as
  required.
  
  Now we assume $C_g = 0$. 
  If $X$ is new this implies $\Dir(X^1,X^2) = \Dir(Y^1,Y^2)$, hence case (ii) can not appear.
  We consider only case (i), $\Dir(Y^1,Y^2) = 0$, implying $C_f=0$ hence $X \leftrightarrow C$. 
  We may assume $Y^1_g = +$ (possibly interchanging superscripts) and $C_p = Y^1_p$ (possibly using $-C$)
  hence $Y^1$ is on the path from $C$ to $Y$.
   If $X$ is new $X^1_g = 0$ would imply $X^1 = C$
   and Proposition \ref{prop:QuadrangleConstellationNoCrossingPoint} would yield 
   $X^1 = Y^1$ contradicting $Y^1 \neq C$. We obtain $X^1 \neq Y^1$,  $X^1_g = +$ and $X^1$  being on the path from $X$ to $C$.
  Let $h \in \supp(Y) \setminus \supp(Y^1)$. 
  
  We assume $X_h = 0$. Then $X \circ Y^1$ is
  an edge by Proposition \ref{prop:edgesBetweenCrossingEdges} and
  directed like $X \circ Y$ by the Triangle Lemma \ref{lem:TriangleLemma}. 
  If $X$ is new Proposition \ref{prop:uniqueEdges2} (i) yields $Y^1 = X^1$ which is impossible.
  Hence $X$ is old in that case and 
  the path $(X,Y^1,Y)$ is as required.    
  Finally, we consider the
  case $X_h = Y_h \neq 0$, recall that $X$ and $Y$ are conformal.
  Cocircuit elimination of $p$ between $-C$ and $Y^1$ yields $Y$ hence
  $-C_h = Y_h = X_h$.  Cocircuit elimination of $h$ between $C$ and
  $X$ yields a cocircuit $Z$ with $Z_g = +$ and $Z_p = C_p = Y^1_p$.
  The path $Q_1$ from $X$ to $Z$ is undirected and starts with $(X,X^1)$ if $X$ is new.
  Proposition \ref{prop:edgesBetweenCrossingEdges} yields that
  $Z$ and  $Y^1$ span a line   and the Zero-line Lemma \ref{lem:ZeroLineLemma} yields that the path $Q_2$ from $Z$ to $Y^1$ is
  directed like the edge $X \circ Y$.  All cocircuits of $Q_1 \setminus X$ and $Q_2$ lie in $(\mathcal{O}',g)$, have the same $p$-value $\neq 0$ hence are
  old. It holds $Y^1 \leftrightarrow Y$ and the path $Q_1 \cup  (Q_2 \setminus Z) \cup Y$ is as required.
\end{proof}

We define normalized paths. 

\begin{definition}
A path / cycle $P$ of cocircuits in $(\mathcal{O}',g)$ is called {\em prenomalized} iff for each pair $(X,Y)$ or $(Y,X)$ of subsequent cocircuits in $P$ 
with $X$ being new, derived from $X^1 \circ X^2$ holds:
\begin{enumerate}
\item If $X^1 \leftrightarrow X^2$ then $Y$ is old and $Y \in \{X^1,X^2 \}$. 
\item If $Y$ is new, derived from $Y^1 \circ Y^2$ and $Y^1_p = X^1_p$ then $\Dir(X^1,X^2) = \Dir(Y^1,Y^2) \neq 0$. 
\end{enumerate}
We call $P$ {\em normalized} if, additionally, $X$ being new, derived from $X^1 \circ X^2$ and $X^1 \leftrightarrow X^2$ implies that $X$ appears only
in triples $(X^i,X,X^{3-i})$ for $i \in \{ 1,2\}$ in $P$.
We call a maximal sequence of new cocircuits in $P$ a {\em new piece of P}.
\end{definition}

We can always derive a directed normalized cycle from a directed cycle.

\begin{lemma}\label{lem:NonZeroDerivedPath}
Let $P = (P^1, \hdots, P^n)$  be a directed path in $(\mathcal{O}',g,f)$.
Then there is a prenormalized path $Q = (P^1, \hdots, P^n)$ in $(\mathcal{O}', g,f)$ directed like $P$ containing all cocircuits of $P$ and additionally only old cocircuits.
If $P$ is a cycle, then $Q$ contains a normalized directed cycle.
\end{lemma}

\begin{proof}
Let the path $P$ have $1 \leq k \leq n-1$ {\em violations} which are pairs of subsequent cocircuits violating the conditions of the definition of a prenormalized path.
Let $(X,Y)$ be a violation in $P$. We have two cases.
\begin{enumerate}
\item Let $X$ be an old cocircuit and $Y$ be new derived from $Y^1 \circ Y^2$ with $Y^1 \leftrightarrow Y^2$ and
$X \notin \{ Y^1,Y^2\}$ (or vice versa). Then Proposition \ref{prop:uniqueEdges2} (i) yields $X_p = 0$ (or $Y_p = 0$) and Lemma \ref{lem:path1lem} 
yields a path $P' = X, \hdots, Y^i,Y$  (or $P' = (X, X^i, \hdots, Y)$) with $i \in \{ 1,2\}$
directed like $X \circ Y$ containing only old cocircuits except $Y$ (or $X$). 
\item Let $X$ and $Y$ be new cocircuits derived from $X^1 \circ X^2$ and $Y^1 \circ Y^2$.
In a violation we necessarily have either $\Dir(X^1,X^2) = 0,\, \Dir(Y^1,Y^2) = 0$ or $\Dir(X^1,X^2) = -\Dir(Y^1,Y^2) \neq 0$.
In all cases Lemma  \ref{lem:path1lem} yields a path $P' = {(X, X^i \hdots, C, \hdots, Y^i, Y)}$ with $i \in \{ 1,2\}$ directed like $X \circ Y$
containing only old cocircuits except $X$ and $Y$.
\end{enumerate}
If we substitute the subsequent cocircuits $X,Y$ in the path $P$ with the path $P'$ we get a path with $k-1$ violations. 
Using induction yields finally the required path $Q$ with no violations.
If $P$ is a cycle then $Q$ is a closed path. For each new cocircuit $X$ in $Q$, 
derived from $X^1 \circ X^2$ with $X^1 \leftrightarrow X^2$
we substitute every triple $(X^i,X,X^i)$ with $i \in \{ 1,2\}$ in $Q$ with $X^i$
and obtain a directed normalized closed path $Q'$. It is easy to see that $Q'$ contains a directed normalized cycle.
\end{proof}

In new pieces of normalized paths there are constraints concerning the f- and g-values.

\begin{prop}\label{cor:dirNormPathNoFgeuqalZero}
  Let $P$ be a normalized path / cycle in $(\mathcal{O}',g,f)$ where
  all cocircuits have the same $f$-value $v_f \neq 0$ and let $X,Y$ be
  two new cocircuits of $P$ such that all cocircuits
    on the path from $X$ to $Y$ are new.  Let $X$ resp.\ $Y$ be
    derived from $X^1 \circ X^2$ and $Y^1 \circ Y^2$ with $X^1_p =
  Y^1_p$.  Then $X^1_f \ne 0$ or $Y^1_g \ne 0$. In particular we never
  have $Y^1_g = Y^1_f = 0$.
\end{prop}

\begin{proof}
The last statement $Y^1_g = Y^1_f = 0$ implies $Y^1 \leftrightarrow Y$ therefore $Y^1$ would be in $P$ because $P$
is normalized contradicting $P$ being in $(\mathcal{O}',g)$.
If we assume $Y^1_g = 0$ (hence $Y^1_f = v_f$) and $X^1_f = 0$ (hence $X^1_g = +$)
we obtain $\Dir(X, X^1) = -v_f = -Y^1_f = -\Dir(Y, Y^1)$ as $Y^1_g=0$ contradicting the fact that $P$ is normalized.
%
\end{proof}

The next proposition follows from the definition of a normalized path / cycle applying Lemma \ref{prop:ChangeOfFAndG1} 
and Corollary \ref{cor:ChangeOfFAndG2}, and from Proposition \ref{prop:DirectedCycleYieldsDirectedCycle}.

\begin{prop}\label{prop:exchangegfstaysnormalized}
Let $P$ be a normalized path / cycle in $(\mathcal{O}',g,f)$ 
where all cocircuits and have the same $f$-value $v_f \neq 0$.
If $v_f = +$ then $P$ is also a normalized path / cycle in $(\mathcal{O}',f,g)$
otherwise $-P = (-P^1, \hdots, -P^n)$.  
\end{prop}

We can reconstruct some directed normalized cycles in other oriented
matroid programs.  This is an application of the Projection
Lemma and of Proposition \ref{prop:projectionLemmaWithUndirectedEdges}.

\begin{lemma}\label{lemma:ProjectionApplication}
Let $P = P^1, \hdots P^n$ be a directed
  normalized cycle in $(\mathcal{O}',g,f)$ with only new cocircuits.
  Then we necessarily have two disjoint sequences $\{ U^1, \hdots, U^n\}$ and
  ${\{V^1, \hdots, V^n\}}$ of cocircuits such that $P^i$ is derived from
  $U^i,V^i$ and $U^i_p = U^j_p = -V^i_p = -V^j_p$ for all $1 \leq i,j
  \leq n$.  If an edge $P^i \circ P^{i+1}$ in $P$ is directed, then let
  $V^i = V^{i+1}$ hold.  Then there is a directed cycle in
  $(\mathcal{O}',p,f)$ or in $(\mathcal{O}',p,g)$.
\end{lemma}

\begin{proof}
By Proposition \ref{cor:dirNormPathNoFgeuqalZero} in $V^1, \hdots, V^n$ we have either $g \neq 0$ or $f \neq 0$ for all cocircuits.
In the first case, cocircuit elimination of $g$ between $-P^i$ and $V^i$ if $V^i_p = +$ or between $-V^i$ and $P^i$ if $V^i_p = -$ yields $W^i$ with $W^i_p = W^j_p = +$.
Then $W = W^1, \hdots, W^n$ is a closed path and is directed in  $(\mathcal{O}',p,f)$ because of the Projection Lemma \ref{lem:ProjectionLemma} and Proposition \ref{prop:projectionLemmaWithUndirectedEdges}. $W$ contains a directed cycle.
In the second case, we know from Proposition \ref{prop:PathGraphWellKnownFacts} (iv) that all cocircuits in $P$ have the same $f$-value $v_f \neq 0$.
If $v_f = +$ let $P' = P$ if $v_f = -$ let $P' = -P$. Then $P'$ is a normalized directed cycle in $(\mathcal{O}',f,g)$ because of the Proposition before.
We eliminate $f$ like in the case before and get a directed closed path hence a directed cycle in $(\mathcal{O}',p,g)$.
\end{proof}

\subsection{The Orientation of the Cocircuit Graph of $(\mathcal{O} \cup p, g,f)$}\label{section:DirGraphgf}

We give a complete description of the orientation of the edges of that graph where $\mathcal{O} \cup p$ is a lexicographic extension.
Lemma \ref{lem:DirPresOfOutsideEdges} and Proposition \ref{prop:derivedDirectedLikeOrig} already handle  the case of edges that are not new.
We look here at the case of new edges. First, we formulate new assumptions.

\begin{definition}
  By {\em Assumptions \textcolor{black}{B}} we fix the following
  situation: Let $\mathcal{O}$ be an oriented matroid of rank $\rk$ with
  groundset $E$ and let $f \in E$.  Let $\textcolor{black}{I} = [e_1,
  \hdots, e_k]$ with $k \leq \rk$ be an ordered set of independent
  elements of $E$, possibly with $f \in \textcolor{black}{I}$.  Let
  $\mathcal{O}' = \mathcal{O} \cup p$ be the lexicographic extension
  $\mathcal{O}[e^+_1, \hdots, e^+_k]$.

Let $X$ be a cocircuit with $X_p = 0$ derived from the two cocircuits $X^1$ and $X^2$.

Let $Y$ be a cocircuit with $Y_p = 0$ derived from the two cocircuits $Y^1$ and $Y^2$.

Let $X_g = Y_g = +$.  Let $X \circ Y$ be an edge and $X,Y$ be conformal cocircuits.
We assume $X^1_p = Y^1_p  = - X^2_p = -Y^2_p \neq 0$.
Let $\ind_X$ and $\ind_Y$ be the indices of $X$ and $Y$ wrt.\ \textcolor{black}{$I$}. 
\end{definition}
\begin{remark}
In this subsection, we omit the case that one (or more) of the cocircuits is zero on $f$ and on $g$.
\end{remark}

\begin{prop}\label{prop:DistinctCocircuitsNotCrossingPoint}
Let Assumptions B hold. If the four cocircuits $X^1,X^2,Y^1,Y^2$ are distinct,
then none of the cocircuits is the crossing point of the edges $X^1 \circ X^2$ and $Y^1 \circ Y^2$.
\end{prop}

\begin{proof}
That is a direct consequence of Proposition \ref{prop:QuadrangleConstellationNoCrossingPoint}.
\end{proof}

\begin{prop}\label{prop:sameIndexEdgesOnTheSamePSide1}
Let Assumptions B hold. 
Let $X$ and $Y$ have the same index $\ind_X = \ind_Y = i$.
Then we have either $X^1_{e_i} = Y^1_{e_i} = 0$ or $X^2_{e_i} = Y^2_{e_i} = 0$.
 \end{prop}

\begin{proof}
  It follows from Lemma \ref{lem:corrCocircuit} that
    for one and only one of the cocircuits $X^1,X^2$ resp,\ $Y^1,Y^2$
    that give rise to $X$ resp.\ $Y$ we have $e_i = 0$.  We assume
  $X^1_{e_i} = Y^2_{e_i} = 0$. Then we have $X^2_p = X^2_{e_i} \neq 0$
  and $Y^1_p = Y^1_{e_i} \neq 0$.  Because $X^2_p = -Y^1_p$ we obtain
  $X^2_{e_i} = -Y^1_{e_i}$ hence $X_{e_i} = -Y_{e_i} \neq 0$ and $X,Y$
  are not conformal, contradicting assumptions B.
\end{proof}

\begin{prop}\label{prop:sameIndexEdgesOnTheSamePSide}
Let Assumptions \textcolor{black}{B} hold. 
Let $i = \ind_X$ and $j = \ind_Y$ with $i < j$.
Then we have $Y^1_{e_i} = Y^2_{e_i} = 0$ and either $X^1 = Y^1$ and $X^2_{e_i} \neq 0$
or $X^2 = Y^2$ and $X^1_{e_i} \neq 0$.
 \end{prop}

\begin{proof}
Let $C$ be the crossing point between the two edges $X^1 \circ X^2$ and $Y^1 \circ Y^2$.
It holds $C_{e_i} = 0$. On the other hand, we have either $X^1_{e_i} = 0$ or  $X^2_{e_i} = 0$.
We may assume $X^1_{e_i} = 0$. Then we have $X^2_{e_i} \neq 0$ and 
$C = X^1$ because of Proposition \ref{prop:unqieCocircuitOnEndge} and $X^1_g \in \{0,+\}$.
But then Proposition \ref{prop:QuadrangleConstellationNoCrossingPoint} yields $Y^1 = X^1$ which shows everything.
\end{proof}

We define all constellations that can appear between two new comodular, conformal cocircuits.

\begin{definition}\label{def:constellations}
Let Assumptions \textcolor{black}{B} hold. 
 
Let $X$ and $Y$ have the same index $\ind_X = \ind_Y = i$. 
We assume $X^2_{e_i} = Y^2_{e_i} = 0$. We have the following constellations.
 \begin{itemize}
\item If $X^1 \neq Y^1$ and $X^2 = Y^2$ we call it {\em constellation A}, 
\item if $X^1 \neq Y^1$ and $X^2 \neq Y^2$ we call it {\em constellation B}.
\item if $X^1 = Y^1$ and $X^2 \neq Y^2$ we call it {\em constellation AB}.
 \end{itemize}
 Let $X$ and $Y$ have different indices. We may asume that $i = \ind_X < \ind_Y = j$ and $X^2_{e_i} = 0$ hence $X^1_{e_i} \neq 0$. Then we have
  $Y^1_{e_i} = Y^2_{e_i} = 0$ and $X^2 = Y^2$.  We call that constellation an {\em index change constellation}.
 If additionally  we have $Y^1_{e_j} = 0$, hence $Y^2_{e_j} \neq 0$ (or vice versa) we call it an an {\em index change constellation with (or without) side change}.
  \end{definition}
 
 \begin{prop}\label{prop:constellationBX2Y2Edge}
Let constellation B or AB hold.
Then $X^2 \circ Y^2$ is an edge and $X^2,Y^2$ are conformal.
\end{prop}

\begin{proof}
Proposition \ref{prop:CuttingConditionInExtension} and Proposition \ref{prop:edgesBetweenCrossingEdges} yield that $X^2 \circ Y^2$ is an edge. 
If $X^2$ and $Y^2$ were not conformal and we would have $X^2_e = +$, and $Y^2_e = -$, then we would also have $X_e = +, Y_e = -$, and therefore $X, Y$ would not be conformal. 
\end{proof}
 
 \begin{definition}
Let Assumptions \textcolor{black}{B} hold and let $i$ be the lower index of the two cocircuits $X$ and $Y$. 
We call $X$, $Y$ the {\em original cocircuits} and $X \circ Y$ the {\em original edge}.
All old cocircuits with $e_i = 0$ in the constellations are corresponding cocircuits  \textcolor{black}{(see Remark \ref{rem:corrCocircuit})}.
If we have two different corresponding cocircuits in a constellation, they form an edge with $e_i = 0$ which we call the {\em corresponding edge} of the constellation.
\end{definition}

\begin{figure}
  \centering
\includegraphics[scale=0.7]{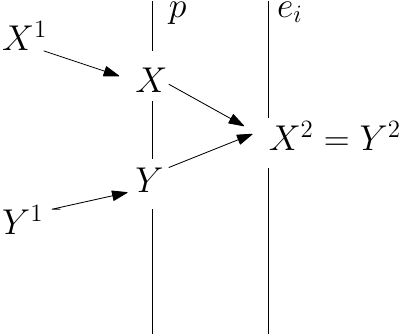} 
  \caption{Constellation A \label{fig:constA}}
\end{figure}

\begin{figure}
  \centering
\includegraphics[scale=0.7]{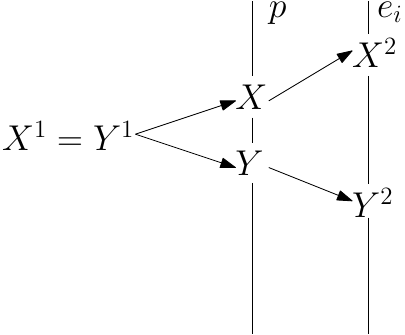} 
  \caption{Constellation AB \label{fig:constAB}}
\end{figure}

\begin{figure}.
  \centering
\includegraphics[scale=0.7]{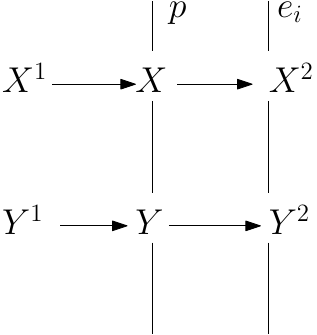} 
  \caption{Constellation B \label{fig:constB}}
\end{figure}

\begin{figure}
  \centering
\includegraphics[scale=0.7]{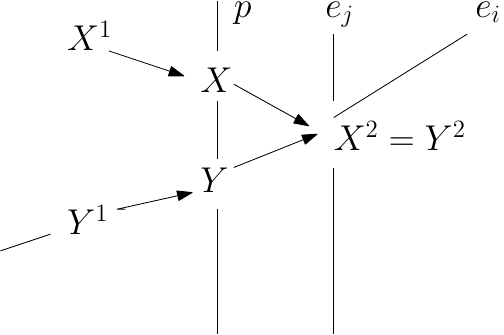} 
  \caption{Index Change Constellation Without Side Change \label{fig:UpconstIndChange}}
\end{figure}

\begin{figure}
  \centering
\includegraphics[scale=0.7]{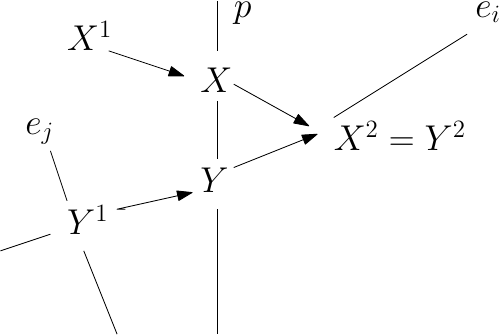} 
  \caption{Index Change Constellation With Side Change \label{fig:DownconstIndChange}}
\end{figure}

We start with a small and clear proposition concerning $g$.

\begin{prop} Let Assumptions \textcolor{black}{B} hold and
let $i$ be the index of a cocircuit with $g = +$. Then $g \notin \{e_1, \hdots, e_{i-1}\}$.
\end{prop}

Now we give a description of the directions of $(X,Y)$ in nearly all constellations. 
We start with the index change constellation.

%
%

\begin{prop}\label{prop:upDownIndexChange}
Let the index change constellation hold. 
We may assume $i = \ind_X  < \ind_Y = j$ and $X^2 = Y^2 = C$ and
let $X_g = Y_g = +$. If $\Dir(C,Y) = 0$, then $\Dir(X,Y) = \Dir(X,C)$ holds. 
Otherwise, we obtain
$\Dir(X,Y) = \Dir(C,Y)$. 
\end{prop}

\begin{proof}
 \textcolor{black}{First we may assume that $C_g = +$.} The first statement follows directly from the Triangle Lemma \ref{lem:TriangleLemma}. Hence let $\Dir(C,Y) \neq 0$. 
If we have $\Dir(C,X) = 0$ or $\Dir(C,X) =  - \Dir(C,Y)$ the Triangle Lemma  \ref{lem:TriangleLemma}  yields $\Dir(X,Y) =  \Dir(C,Y)$.
Hence we may assume $\textcolor{black}{0 \neq } \Dir(C,X) = \Dir(C,Y)$.
We have $C_{e_i} = Y_{e_i} = 0$ and $X_{e_i} \neq 0$.
Cocircuit elimination of $g$ between $-X$ and $C$ yields
the cocircuit $X^3$ and between $-Y$ and $C$ yields the cocircuit $Y^3$. We have $Y^3_f \neq 0$.
We have also $i = \ind_{X^3} < j \textcolor{black}{ \leq }  \ind_{Y^3}$.
The Projection Lemma \ref{lem:ProjectionLemma} yields  $\Dir_{g,f}(X,Y) = -\Dir_{p,f}(X^3,Y^3)$.
Lemma \ref{lem:case6lexExt} yields $\Dir_{p,f}(X^3,Y^3) = Y^3_f \neq 0$. But it holds also $Y^3_f = \Dir_{g,f}(Y,C)$.
We obtain $\Dir_{g,f}(X,Y) = - \Dir_{g,f}(Y,C)$. 
 \textcolor{black}{Now we assume that $C_g = 0$ hence $C_f = v_f \neq 0$. We obtain directly $X_f = Y_f = v_f$ and get 
 the statement of the Proposition by exchanging $g$ and $f$ and using Lemma \ref{prop:ChangeOfFAndG1} or Corollary \ref{cor:ChangeOfFAndG2}.}
\end{proof}

%
%

The next lemma handles constellation B and AB.
We use the Projection Lemma \ref{lem:ProjectionLemma} quite often in the proof.

\begin{lemma}\label{lem:constellationBDirCorrEdge}
Let constellation B or AB hold. Let $X_g=Y_g=+$ and 
let $i$ be the index of the cocircuits $X$ and $Y$. Let $f,g \notin \{e_1, \hdots, e_{i-1}\}$. 
Let $X^2_p = Y^2_p = \textcolor{black}{v}_p \neq 0$ and let $\Dir(X,X^2) = \Dir(Y,Y^2) \neq 0$.
\textcolor{black}{Let $X^2_g \neq 0$ or $Y^2_g \neq 0$ or both and}
let $j$ be the index of the cocircuit $Z$ yielded by cocircuit elimination of $g$ between $-X^2$ and $Y^2$. Then the following holds:
\begin{enumerate}
\item  If $\Dir(X^2,Y^2) \neq 0$ then $\Dir(X,Y) = \Dir(X^2,Y^2)$.  
\item   If $\Dir(X^2,Y^2) = 0$ and $j \leq k$ then \[\Dir(X,Y) = - \textcolor{black}{v}_p * \Dir(X,X^2) * \Dir_{g,e_j}(X^2,Y^2). \]
\item  If $\Dir(X^2,Y^2) = 0$ and $j = k+1$ then $\Dir(X,Y) = 0$.
\end{enumerate}
\textcolor{black}{Now let $X^2_g = Y^2_g = 0$ and let $j$ be the index of the cocircuit $Z$ yielded by cocircuit elimination of $f$ between $-X^2$ and $Y^2$. Then the following holds:
\begin{enumerate}
 \setcounter{enumi}{3}
\item    If $j \leq k$ then \[\Dir(X,Y) = - v_p * \Dir(X,X^2) * \Dir_{f,e_j}(X^2,Y^2). \]
\item  If $j = k+1$ then $\Dir(X,Y) = 0$.
\end{enumerate}}
If $f$ \textcolor{black}{and $g$} are in general position and not in $\textcolor{black}{I}$ only case $\textcolor{black}{(i)}$ can appear.
\end{lemma}
\begin{proof}
  Let $\Dir(X,X^2) = v_f$. We have $X_f = Y_f \neq 0$ because $v_f \neq
  0$ and $j > i$ because $X^2$ and $Y^2$ have index $> i$. We assume
$X^2_g = Y^2_g = +$. Cocircuit elimination of $g$ between $-X$ and
$X^2$ yields $X^3$. Cocircuit elimination of $g$ between $-Y$ and
$Y^2$ yields $Y^3$. We have $X^3_f = Y^3_f = v_f \neq 0$.
Let $C$ be the crossing point between the two edges $X
  \circ X^2$ and $Y \circ Y^2$. We assume $C_g = +$ hence $C \neq X^3
  \neq Y^3 \neq C$ and postpone the other case at the end.  We have
  $v_p = X^3_p = Y^3_p= X^3_{e_i} =Y^3_{e_i} \neq 0$ because $X^3$ and
  $Y^3$ have index $i$.  Let $\epsilon = -$ if $C_p = X^2_p$ and
  $\epsilon = +$ otherwise. The Projection Lemma
\ref{lem:ProjectionLemma} yields
\begin{align}
&\Dir_{g,f}(X,Y) = \epsilon * \Dir_{p,f}(X^3,Y^3)  \label{eq:aligned1} \\
&\Dir_{g,f}(X^2,Y^2) = \epsilon * \Dir_{e_i,f}(X^3,Y^3) \text{ if } f \neq e_i \label{eq:aligned2} \\
&\Dir_{g,e_j}(X^2,Y^2) = \epsilon * \Dir_{e_i,e_j}(X^3,Y^3)  \label{eq:aligned3} 
\end{align}
If $f \neq e_i$ and $\Dir_{g,f}(X^2,Y^2) \neq 0$ equation \ref{eq:aligned2}  yields $\Dir_{e_i,f}(X^3,Y^3) \neq 0$. 
Using equation \ref{eq:aligned1} , Lemma \ref{lem:case3lexext} and equation \ref{eq:aligned2}  we have 
\[ \Dir_{g,f}(X,Y) = \epsilon * \Dir_{p,f}(X^3,Y^3)  = \epsilon * \Dir_{e_i,f}(X^3,Y^3) =  \Dir_{g,f}(X^2,Y^2) \]
Otherwise the assumptions of Lemma \ref{lem:case4lexExt}  hold for $X^3$ and $Y^3$. We assume $i < j \le k$. Then $j$ is also the smallest index with $\Dir_{e_i,e_j}(X^3,Y^3) \neq 0$ because equation \ref{eq:aligned3} holds for all $j > i$. We obtain (using equation \ref{eq:aligned1}, Lemma \ref{prop:ChangeOfFAndG1} or Corollary \ref{cor:ChangeOfFAndG2}, Lemma \ref{lem:case4lexExt} and equation \ref{eq:aligned3}):
\begin{align*}
\Dir_{g,f}(X,Y) &= \epsilon * \Dir_{p,f}(X^3,Y^3) = - \epsilon * v_p * X^3_f * \Dir_{f,p}(X^3,Y^3) \\ 
&=  - \epsilon * v_p * X^3_f * \Dir_{e_i,e_j}(X^3,Y^3) = - v_p * X^3_f * \Dir_{g,e_j}(X^2,Y^2)
\end{align*}
If $j = k+1$ we have 
\[  \Dir_{g,f}(X,Y) = \epsilon * \Dir_{p,f}(X^3,Y^3) = 0.  \]
Now let $X^2_g = +$ and $Y^2_g = 0$. 
The Projection Lemma \ref{lem:ProjectionLemma} and Lemma \ref{lem:case6lexExt} yield
 \[\Dir_{g,f}(X,Y) = \Dir_{p,f}(X^3,Y^2) = Y^2_f = \Dir_{g,f}(X^2,Y^2). \]
The case $X^2_g = 0$ and $Y^2_g = +$ works analogously. 
If $X^2_g = Y^2_g = 0$ we have $X^2_f = Y^2_f = \Dir(X,X^2) = v_f$. The
Projection Lemma \ref{lem:ProjectionLemma}, Lemma
\ref{prop:ChangeOfFAndG1} and the fact that $Z$ is old (Lemma
\ref{lem:DirPresOfOutsideEdges}) having index $j$ yield
\[ \Dir_{g,f}(X,Y) = \Dir_{p,f}(X^2,Y^2) =  - v_p * v_f * \Dir_{f,p}(X^2,Y^2) =  - v_p * v_f * Z_p \]
This equals $0$ if $j = k+1$ otherwise we obtain 
\[ \Dir_{g,f}(X,Y)  = - v_p * v_f * Z_{e_j}  = - v_p * v_f * \Dir_{f,e_j}(X^2,Y^2).\]
Finally if $X^3 = Y^3= C$ it must be $X^3_f \neq 0$.
 Because we have $X_f = Y_f \neq 0$ we can exchange $f$ and $g$ using Lemma \ref{prop:ChangeOfFAndG1} or Corollary \ref{cor:ChangeOfFAndG2} and work like before.
 In case (ii) we need Proposition  \ref{prop:TransitivityOfDir} yielding $\Dir_{g,e_j}(X^2,Y^2) =  \Dir_{f,e_j}(X^2,Y^2)$ if $\Dir(X^2,Y^2)  = 0$.
 The last statement is obvious.
 \end{proof}

To be complete, we add a short proposition concerning constellation A which follows directly from the  Projection Lemma \ref{lem:ProjectionLemma} and from Proposition \ref{prop:ChangeOfFAndG1}.

\begin{prop}\label{prop:gConstellationABB}
Let constellation A hold. Let $X_f = Y_f \neq 0$ and $X_g = Y_g = +$.
Let $X^2_g = +$ then cocircuit elimination of $g$ between $-X$ ($-Y$) and $C$ yields $X^3$ ($Y^3$). 
We obtain $\Dir_{g,f}(X,Y) = - \Dir_{p,f}(X^3,Y^3)$. If $X^2_g = 0$ and $X^2_f = X_f$ we eliminate $f$ 
between $-X$ ($-Y$) and $C$ obtaining $X^3$ ($Y^3$).  We have $\Dir_{g,f}(X,Y) =  - v_f * \Dir_{f,g}(X,Y) =  v_f * \Dir_{p,g}(X^3,Y^3)$.
\end{prop}


\subsection{The Corresponding Cycle}

Now we construct to a normalized cycle $P$ its corresponding cycle $Q$
containing only old, corresponding cocircuits (recall
Remark~\ref{rem:corrCocircuit}). We show that in most cases $Q$ is
directed in a program of $\mathcal{O}$ iff $P$ is.  In this chapter,
let $\mathcal{O}$ be an oriented matroid of rank $\rk$ with groundset
$E$ and $f \in E$.  Let $\textcolor{black}{I} = [e_1, \hdots, e_k]$
with $k \leq \rk$ be an ordered set of independent elements of $E$.
Let $\mathcal{O}' = \mathcal{O} \cup p$ be the lexicographic extension
$\mathcal{O}[e^+_1, \hdots, e^+_k]$.  We construct the {\em
  corresponding path}.

\begin{definition}\label{pathConstellation222}
Let $P = (P^1, \hdots, P^n)$ be a path / cycle in $\mathcal{O}'$.
Let $Q^i$ for $1 \leq i \leq n$ such that $Q^i = P^i$ if $P^i$ is an old cocircuit
and $Q^i = X$ if $P^i$ is a new cocircuit and $X$ its corresponding cocircuit.
We call $Q' = (Q^1, \hdots, Q^n)$ the {\em corresponding sequence} to $P$. In $Q'$ we may have subsequent repeating cocircuits.
Identifying subsequent equal cocircuits, we call $Q$ the {\em corresponding path} to $P$. 
\end{definition}

We collect some facts about the corresponding sequence / path.
\textcolor{black}{In the following if we say an edge is directed, we mean directed in $(\mathcal{O}',g,f)$ or $(\mathcal{O},g,f)$.}

\begin{prop}\label{prop:subsequentOriginalCocircuits}
Let $P$ be a path in $(\mathcal{O}',g)$ and let 
\textcolor{black}{$Q'$} be the corresponding \textcolor{black}{sequence} to $P$. Let $P^i, P^{i+1}$ be two subsequent cocircuits in $P$.
Then \textcolor{black}{ we have $Q^i \neq Q^{i+1}$ in $Q'$ in four cases.}
\begin{enumerate}
\item $P^i, P^{i+1}$ are both old. 
\item One of both cocircuits is new, the other one is old but not corresponding. 
\item $P^i, P^{i+1}$ are in constellation $AB$ or $B$.
\item $P^i, P^{i+1}$ are in an index change constellation with side change. 
\end{enumerate}
 \textcolor{black}{We omit the case that $Q^i_g = Q^i_f = 0$ or $Q^{i+1}_g = Q^{i+1}_f = 0$ or both.}
 \textcolor{black}{$F = Q^i \circ Q^{i+1}$ and $E = P^i \circ P^{i+1}$ are edges.}
\textcolor{black}{In the four cases, we obtain}
\begin{enumerate}
\item $F$ is directed like $E$.
\end{enumerate}
\begin{enumerate}
\item [(ii)] $F$ is directed like $E$.
\item [(iii)]  $F$ is undirected or directed like $E$.
\item [(iv)] $F$ is undirected or directed like $E$. If $P$ is normalized, $F$ is directed like $E$.
\end{enumerate}
\end{prop}

\begin{proof}
If $P^i,P^{i+1}$ are old we have $E = F$, if only one cocircuit $X \in \{P^i,P^{i+1} \}$ is new, then $F$ is the edge where 
$X$ is derived from, and Proposition \ref{prop:derivedDirectedLikeOrig} yields (ii).
Lemma \ref{lem:constellationBDirCorrEdge} yields (iii) and Proposition \ref{prop:upDownIndexChange} (iv). In all cases $F$ is an edge, this is Proposition \ref{prop:constellationBX2Y2Edge} and Theorem \ref{theorem:ext}.
\end{proof}

Inversely, we have.

\begin{prop}\label{prop:subsequentCorrespondingCocircuits}
Let $P$ be a path in $(\mathcal{O}',g)$ and let 
$Q$ be the corresponding path to $P$. Let $Q^j, Q^{j+1}$ be two subsequent cocircuits in $Q$.
Then they correspond to two subsequent cocircuits $P^i$ and $P^{i+1}$ in $P$.
The cocircuits  $Q^j \setminus p, Q^{j+1} \setminus p$ are conformal in $\mathcal{O}$ and $F = Q^j \setminus p \circ Q^{j+1} \setminus p$ is an edge in $\mathcal{O}$.
\textcolor{black}{We call $E = P^i \circ P^{i+1}$  {\em the original edge} and $F$ the {\em corresponding edge} to $E$.}
If we assume $Q^j_g = Q^{j+1}_g = +$ and if $F$ is directed then it is directed like $E$. 
\end{prop}

\begin{proof}
The existence of $P^i$ and $P^{i+1}$ follows from the definition.
If $Q^j \setminus p, Q^{j+1} \setminus p$ were not conformal, then also not $P^i$ and $P^{i+1}$.
The other statements follow from the Proposition before.
\end{proof}

We prove some facts about the corresponding path.

\begin{prop}\label{prop:propsCorresPath}
Let $P = (P^1, \hdots P^n)$ be a path in $\mathcal{O}'$
and $Q$ its corresponding path. The following holds:
\begin{enumerate}
\item $Q$ contains only old cocircuits and contains all old cocircuits of $P$.
\item If $P^i$ is a new cocircuit, then its corresponding cocircuit lies in $Q$.  

\item $Q$ is a path in $\mathcal{O}$. If $P$ is a cycle, then $Q$ is a closed path.  
\item If an edge in $Q$ is directed, the original edge in $P$ has the same direction.
\item If $P$ is directed, then $Q$ is directed like $P$ or undirected.
\end{enumerate}
\end{prop}

\begin{proof}
(i), (ii) and (iii) follow from the construction of $Q$, (iv) is the Proposition before and
(v) follows from (iv).
\end{proof}

Now we consider normalized paths.
\textcolor{black}{Recall Proposition \ref{cor:dirNormPathNoFgeuqalZero} that in normalized paths none of the corresponding cocircuits has $f = g = 0$.}


%

\begin{lemma}\label{lemma:AllfOrNobody}
Let $P = (P^1, \hdots P^n)$ be a directed normalized closed path in $\mathcal{O}'$
and $Q$ its corresponding closed path. Then
either all cocircuits of $Q$ have $f = 0$ or none of them and
either all cocircuits of $Q$ have $g = 0$ or none of them.
\end{lemma}

\begin{proof}
Because $P$ is directed and because of Proposition  \ref{prop:PathGraphWellKnownFacts} (iv),  all cocircuits in $P$ have the same f-value $v_f$.
Hence, all cocircuits in $Q$ have $f = 0$ or $f = v_f$ and $g = 0$ or $g = +$ (but never $f = g = 0$ because $P$ is normalized).
First, we assume having cocircuits in $Q$ with $g= 0$ (hence $f = v_f$) and other cocircuits with $g= +$.
But then we have at least one edge going from $g = 0$ to $g = +$ having direction $-v_f$ and one edge going from $g = +$ to $g = 0$ with direction $v_f$
(see the definition of the Dir-function), hence also $P$ would have edges with opposite directions because of Proposition \ref{prop:propsCorresPath}.
If we assume all cocircuits having $g=0$ then all cocircuits have $f= v_f$ and we are done.
Lastly, we assume all cocircuits having $g=+$. If we had a cocircuit with $v_f$ and a cocircuit with 
$0$ in $Q$, we would get at least two edges with opposite directions and hence also in $P$  like before which is impossible. 
\end{proof}

\begin{prop}\label{prop:CorrPathToCycleWithOldCocNotDirected}
The corresponding closed path $Q$ to a directed normalized closed path $P$ with old cocircuits has no cocircuits with $f=0$ and with $g=0$.
\end{prop}

\begin{proof}
All cocircuits in $P$ have the same $f$-value $v_f \neq 0$.
\textcolor{black}{Because $P$ contains old cocircuits $Q$ contains also at least one cocircuit with $g=+$ and $f= v_f$ hence all
cocircuits in $Q$ must have $g=+$ and $f=v_f$ because of the Lemma before.}
%
%
 \end{proof}

%

\textcolor{black}{Now we can handle the index change constellations with side change.}

\begin{lemma}\label{lemma:CorrespondingCycleIndexChangesDirected}
The corresponding cycle to a directed normalized cycle with index change constellations with side change is always directed.
\end{lemma}

\begin{proof}
  Let $F$ be a corresponding edge to an index change constellation
  with side change.  Then at least one of the two cocircuits of the
  edge must have $f \neq 0$ and one $g = +$ hence Proposition
  \ref{prop:CorrPathToCycleWithOldCocNotDirected} and Lemma
  \ref{lemma:AllfOrNobody} yield that all $f \neq 0$ and $g = +$ in
  the corresponding cycle.  Proposition
  \ref{prop:subsequentOriginalCocircuits} (iv) yields that $F$ is
  directed like the original edge, hence Proposition
  \ref{prop:propsCorresPath} yields that the corresponding cycle is
  directed.
\end{proof}

We show a first main result concerning the corresponding cycle.

\begin{lemma}\label{lemma:NormalizedWithOutsideCocircuits}
The corresponding closed path $Q$ to a directed normalized closed path $P$ with old cocircuits is always directed.
\end{lemma}

\begin{proof}
If $P$ contains only old cocircuits, then $Q= P$.
Hence, we may assume that there are also new cocircuits in $P$. All cocircuits in $P$ have the same $f$-value $v_f \neq 0$.
Proposition \ref{prop:CorrPathToCycleWithOldCocNotDirected} 
yields that all cocircuits in $Q$ have $g = +$ and $f = v_f$.
Because of Lemma \ref{lemma:CorrespondingCycleIndexChangesDirected} 
we may assume that we have no index change constellations with side change in the cycle.

We show that $Q$ contains at least one directed edge (then it is directed because of Proposition \ref{prop:propsCorresPath} (vi)).
We may assume that $P$ hence $Q$ starts with an old cocircuit.
We consider the first transition in $P$ from an outside cocircuit $P^{i-1}$ to an inside cocircuit $P^i$.
It comes now a pathpiece with new cocircuits derived from cocircuits of two partitions $X = \{X^i \hdots\} $ and $Y = \{Y^i \hdots \}$ with $X^i_p = -Y^j_p$.  
The corresponding path of that inside pathpiece contains only the edges of one partition. We may assume of partition $Y$.
 If $P^{i-1} = X^i$ we have an edge from $X^i$ to $Y^i$ in the corresponding cycle 
that is directed, too. Hence, we may assume that is $P^{i-1} = Y^i$. Because $P$ is normalized, we have $\Dir(Y^i,P^i) = v_y$ with $v_y \in \{+,-\}$.
But then let $P^{i+k}$ be the last new cocircuit in that pathpiece and $P^{i+k+1}$ the next outside cocircuit.
Because $P$ is normalized, we have $\Dir(P^{i+k},Y^{i+k}) = -v_y$ hence because $P$ is directed, $P^{i+k+1}$
must be the cocircuit $X^{i+k}$ which means that we have the directed edge $Y^{i+k} \circ X^{i+k}$ in the corresponding path $Q$.
\end{proof}

\begin{lemma}\label{lem:AllCocircuitsSameIndex}
Let $P = (P^1, \hdots, P^n = P^1)$ be a normalized directed cycle of new cocircuits without index change constellations with side change.
Then all cocircuits in $P$ have the same index.
\end{lemma}

\begin{proof}
Assume to the contrary that in the cycle we have at least one index-change constellation without side change from a cocircuit with lower index to 
a cocircuit with higher index and vice versa.
Proposition \ref{prop:upDownIndexChange} yields then two edges in $P$ in opposite direction because $P$ is normalized.
\end{proof}

\begin{lemma}\label{lemma:DrectedCycleBecauseConstellABandB}
Let $P = (P^1, \hdots, P^n = P^1)$ be a normalized directed cycle of new cocircuits without index change constellations with side change such that at least one edge in constellation B or AB is directed. Then there is a directed cycle in $(\mathcal{O},g,f)$ or in a $(\mathcal{O},g,e_j)$ (or in $(\mathcal{O},f,g)$ or in a $(\mathcal{O},f,e_j)$) where $j$ is higher than the index of  the cocircuits of $P$. \textcolor{black}{The second case appears only if $f$ (or $g$) is not in general position.}
\end{lemma}

\begin{proof}
All cocircuits of $P$ have the same f-value $v_f \neq 0$. Because of Lemma \ref{lem:AllCocircuitsSameIndex} all cocircuits in $P$ have the same index $i$.
We assume $g = +$ for all cocircuits in $Q$. 
\textcolor{black}{Proposition \ref{prop:propsCorresPath} yields that if $Q$ has at least one directed edge in $(\mathcal{O},g,f)$ it is a directed cycle there.}
If $Q$ is undirected then let $j > i$ be the smallest index reached by cocircuit elimination of $g$ between the corresponding cocircuits $-Q^1, Q^2$ of 
the directed edges in constellation B or AB. 
Because of Lemma \ref{lem:constellationBDirCorrEdge} \textcolor{black}{(ii) it must hold $j < r+1$, because $P$ is directed} and $Q$ is directed in $(\mathcal{O},g,e_j)$ because $P$ is normalized and the directions $\Dir(X,X^2)$ are all the same.
\textcolor{black}{If we have $g=0$ for all cocircuits in $Q$ we obtain $f = v_f$ for all of them (because $P$ is normalized). Then $P$ (or $-P$ if $v_f = -$) is also a directed and normalized cycle in $(\mathcal{O},f,g)$, see Proposition \ref{prop:exchangegfstaysnormalized},
and we obtain a directed cycle in $(\mathcal{O},f,g)$ or in $(\mathcal{O},f,e_k)$ with $k > i$ like before.}
\end{proof}

\begin{lemma}\label{lemma:LastCycleProjection}
Let $(P^1, \hdots, P^n = P^1)$ be a normalized directed cycle of new cocircuits without index change constellations with side change such that none of the edges in constellations AB  and B are directed.
Then there is a directed cycle in $(\mathcal{O}',p,f)$ or in $(\mathcal{O}',p,g)$.
\end{lemma}

\begin{proof}
That is Lemma \ref{lemma:ProjectionApplication}. If there is $f=0$ for all cocircuits, we project to $g$ and if $g=0$ we project to $f$.  
\end{proof}

\subsection{Second Main Theorem}\label{subsection:secondMainResult}

We put everything together for the second main Theorem.

\begin{theorem}\label{Theo:SecondMaintheorem2}
Let $\mathcal{O}$ be an  oriented matroid of rank $\rk$ with groundset $E$ and $f \neq g \in E$. 
Let $\textcolor{black}{I} = [e_1, \hdots, e_k]$ with $k \leq \rk$ be an ordered set of independent elements of $E$.
Let $l  \leq k$ be the smallest index
such that $f \in \cl(\{e_1, \hdots, e_l\})$. If there is no such $l$ then let $l = k+1$.
Let $m \leq k$ be the smallest index
such that $g \in \cl(\{e_1, \hdots, e_m\})$. If there is no such $m$ then let $m = k+1$.
\begin{enumerate}
\item Let  $(\mathcal{O},g,f)$ be a Euclidean oriented matroid program and
\item let  $(\mathcal{O},e_1,f), \hdots, \textcolor{black}{(\mathcal{O},e_k,f)}$ be Euclidean oriented matroid programs and
\item let  $(\mathcal{O},e_1,g), \hdots,  \textcolor{black}{(\mathcal{O},e_k,g)}$ be Euclidean oriented matroid programs and
\item if $min(l,m) < k+1$ or if $f$ or $g$ are not in general position let $n =min(max(l,m),k-1)$ and 
let $(\mathcal{O} / \{e_1, \hdots, e_{i-1}\}, e_i, e_j)$
be Euclidean oriented matroid programs for all $1 \le i \le n$ and $i < j \le k$. 
\end{enumerate}
Let $\mathcal{O}' = \mathcal{O} \cup p$ be the lexicographic extension $\mathcal{O}[e^+_1, \hdots, e^+_k]$.  
Then $(\mathcal{O}',g,f)$ is a Euclidean oriented matroid program.
\end{theorem}

\begin{proof}
  First, if we have a directed cycle in $(\mathcal{O}',g,f)$ with old
  cocircuits, we also have a directed normalized closed path in
  $(\mathcal{O}',g,f)$ with old cocircuits because of Lemma
  \ref{lem:NonZeroDerivedPath}. But because of Lemma
  \ref{lemma:NormalizedWithOutsideCocircuits}, we obtain a directed
  corresponding closed path and Proposition
  \ref{prop:PathGraphWellKnownFacts} (i) yields a directed cycle in
  $(\mathcal{O},g,f)$, which is impossible.  Hence, we assume a directed
  cycle $P$ with only new cocircuits. Again, Lemma
  \ref{lem:NonZeroDerivedPath} yields a normalized directed closed
  path. Again, that closed path cannot have old cocircuits. Hence, we
  can assume that $P$ is a directed normalized cycle with only new
  cocircuits.  Let $Q$ be its corresponding closed path.

Because of Proposition \ref{prop:PathGraphWellKnownFacts} (iv) all cocircuits in $P$ have the same $f$-value $v_f$.
We know from Lemma \ref{lemma:AllfOrNobody} that for all cocircuits in $Q$ we have $g=0$ or for none of them.
We assume $g=+$ for all of them; otherwise we have $f=v_f$ for all cocircuits and we look at the cycle $P$ (or $-P$ if $v_f = -$) in $(\mathcal{O}',f,g)$. It stays a directed
normalized cycle there \textcolor{black}{because of Proposition \ref{prop:exchangegfstaysnormalized}.}
We know from Lemma \ref{lemma:CorrespondingCycleIndexChangesDirected}
 that we do not have index change constellations with side change in $P$.
 Hence, from Lemma \ref{lem:AllCocircuitsSameIndex} we know that all cocircuits in $P$ have the same index.
 If in $P$ at least one of the constellation AB or B edges are
 directed, the closed path $Q$ is also directed in $(\mathcal{O},g,f)$
 or in $(\mathcal{O},g,e_j)$ where $j$ is higher than the index of the
 cocircuits in $P$, see Lemma
 \ref{lemma:DrectedCycleBecauseConstellABandB}.  If in $P$ the
 constellation AB and B edges are all undirected \textcolor{black}{(or
   there are no such edges)} and an edge of constellation A is
 directed, the projected closed path $Q'$ is also directed in
 $(\mathcal{O}',p,f)$ (or in $(\mathcal{O}',p,g)$ if we exchange $g$
 and $f$), see Lemma \ref{lemma:LastCycleProjection},
contradicting that $(\mathcal{O}',p,f)$ and
   $(\mathcal{O}',p,g)$ are Euclidean oriented matroid programs,
   as shown by the first main Theorem \ref{theo:lexExtStaysEucl2}.
\end{proof}

\begin{remark}
We use assumption (iv) of the Theorem only to recur on Theorem \ref{theo:lexExtStaysEucl2} hence we can always substitute $e_i$ by $f$ or by $g$ there (see 
Remark \ref{rem:SubstitutionMainTheoremAssumptions}).
\end{remark}

We formulate the Theorem if $g$ and $f$ are in general position.

\begin{corollary}
Let $\mathcal{O}$ be an oriented matroid of rank $\rk$ with groundset $E$. 
Let $\textcolor{black}{I} = [e_1, \hdots, e_k]$ with $k \leq \rk$ be an ordered set of independent elements of $E$.
Let $f \neq g  \in E \setminus \textcolor{black}{I}$ be in general position.
\begin{enumerate}
\item Let  $(\mathcal{O},g,f)$ be a Euclidean oriented matroid program and
\item Let  $(\mathcal{O},e_1,f), \hdots, (\mathcal{O} / \{e_1, \hdots, e_{k-1}\},e_k,f)$ be Euclidean oriented matroid programs and
\item Let  $(\mathcal{O},e_1,g), \hdots, (\mathcal{O} / \{e_1, \hdots, e_{k-1}\},e_k,g)$ be Euclidean oriented matroid programs.
\end{enumerate}
Let $\mathcal{O}' = \mathcal{O} \cup p$ be the lexicographic extension $\mathcal{O}[e^+_1, \hdots, e^+_k]$.  
Then $(\mathcal{O}',g,f)$ is a Euclidean oriented matroid program.
\end{corollary}

\begin{proof}
That follows directly from \textcolor{black}{Theorem \ref{Theo:SecondMaintheorem2}
and from the fact that we needed assumptions (ii) and (iii) in the proof of that theorem only where we applied Corollary \ref{cor:lexExtStaysEucl2} (at the end of the proof) and where we applied the second case of Lemma \ref{lemma:DrectedCycleBecauseConstellABandB}. That case does not appear if $f$ and $g$ are in general position.
Hence, we can replace these assumptions with the weaker assumptions of Corollary \ref{cor:lexExtStaysEucl2}.}
\end{proof}

From that follows directly:

\begin{thmn}[\ref{theorem:SecondMainTheorem}]
A lexicographic extension of a Euclidean oriented matroid
stays Euclidean.
\end{thmn}

 \section{Concluding Remarks}
 
With our result, we can show that if an oriented matroid is Euclidean, it is also Mandel.
Because of Theorem \ref{theorem:SecondMainTheorem} any lexicographic extension fulfills the desired assumptions. 
 The inverse direction is mentioned to be false in \cite{Knauer} but we did not find an explicit example of an oriented matroid
 that is Mandel but not Euclidean in the literature. In a forthcoming paper, we construct such an example explicitly and  show some properties about mutations in Euclidean oriented matroids. 
 
 In the case of uniform oriented matroids, it is easy to see (in fact, this is Lemma \ref{lem:case3lexext} together with Theorem  \ref{theorem:EuclideanessStays}, see also \cite{SturmfelsZiegler}, Proposition 4.4 and Proposition 4.7) that for an inseparable pair $(f,f')$ of elements of the groundset $(\mathcal{O},g,f')$ remains Euclidean if $(\mathcal{O},g,f)$ is. It would be interesting to prove this fact in general (if possible) and to develop it further to ask whether the construction of the {\em connected sum} introduced by Richter-Gebert in \cite{Richter-Gebert}, Chapter 1, preserves Euclideaness. That construction uses {\em flagged} oriented matroids (see \cite{Richter-Gebert}, Definition 1.1), which in a way generalizes the concept of inseparable elements. It is used in \cite{Richter-Gebert}, chapter 2, to build the $R(20)$, an oriented matroid having a mutation-free element, hence not being Mandel (see \cite{Knauer}). In a forthcoming paper, we will use inseparability and techniques similar to here to show that each oriented matroid program of the $R(20)$ has Euclidean and Non-Euclidean regions (supercells). We remark here that Theorem \ref{theorem:SecondMainTheorem} holds 
 also for Euclidean regions in (maybe Non-Euclidean) oriented matroids.

 In a Non-Euclidean oriented matroid program the cocircuit graph has always at least one {\em very strong component}, which means all cocircuits of the component (which must be more than one) are connected to each other in both ways in the graph.
It could be interesting to find further characterizations of these components via the fact that each directed cycle in a component can be normalized, which is Lemma \ref{lem:NonZeroDerivedPath} of this paper.
 
 We are not aware of any result concerning the portions of Euclidean Oriented Matroids which are linear or the asymptotics of the quotient of the number of non-linear versus linear matroids at all. Peter Nelson gave a nice asymptotic result for ordinary matroids~\cite{nelson}. Therefore, we ask.

 \begin{problem}
     If $L_n$ denotes the number of linear oriented matroids on $n$ elements, $E_n$ the number of Euclidean oriented matroids and $O_n$ the number of oriented matroids on $n$ elements. Determine the asymptotics
     \[\frac{L_n}{E_n},\,\frac{E_n}{O_n},\,\frac{L_n}{O_n} \text{ as }n \to \infty\]
 \end{problem}
 
\bibliographystyle{plain} 

\begin{thebibliography}{99}
\bibitem{1} \textsc{W.\ Hochst\"attler}, \textit{Topological Sweeping in Oriented Matroids.} Technical report (2016) 
\bibitem{2} \textsc{A.\ Bj\"orner, M.\ Las Vergnas, B.\ Sturmfels, N.\ White, and G.M.\ Ziegler} \\ \textit{Oriented Matroids (2nd ed.), Encyclopedia of Mathematics and Its Applications, Vol. 46, Cambridge University Press, Cambridge, 1999}
\bibitem{3} \textsc{R.\ Cordovil and K.\ Fukuda}  \textit{Oriented Matroids and Combinatorial Manifolds, European Journal of Combinatorics, 1-14, 1993} 
\bibitem{5} \textsc{A.\ Mandel}  \textit{Topology of Oriented Matroids, Ph D Thesis of A. Mandel, University of Waterloo, 1982} 
\bibitem{6} \textsc{K.\ Fukuda}  \textit{Oriented Matroid Programming, Ph D Thesis of K. Fukuda, University of Waterloo, 1982} 
\bibitem{8} \textsc{H.\ Bruggesser, P.\ Mani}  \textit{Shellable Decompositions of Cells and Spheres., MATHEMATICA SCANDINAVICA, 29, 197-205. , 1971} 
\bibitem{Blandpivot} \textsc{R.G.\ Bland} \textit{New Finite Pivoting Rules for the Simplex Method, Mathematics of Operations Research Vol.\ 2 No.\ 2, pp 103 -- 107}
\bibitem{McMullen} \textsc{P.\ McMullen} \textit{The number of faces of simplicial polytopes, Israel Journal of Mathematics 9 Vol. 4, 559 -- 570, 1971}
\bibitem{Athanasiadis} \textsc{C.\ Athanasiadis} \textit{Zonotopal Subdivisions of Cyclic Zonotopes, Geometriae Dedicata Vol. 86, no 1-3, pp 37 - 57, 2001}
\bibitem{Richter-Gebert} \textsc{J.\ Richter-Gebert} \textit{Oriented Matroids With Few Mutations, Discrete Comp.  Geom. 10, pp 251 - 269, 1993}
\bibitem{Knauer} \textsc{K.\ Knauer and T.\ Marc} \textit{Corners and simpliciality in oriented matroids and partial cubes, European Journal of Combinatorics
Vol.\ 112, 2023, 103714}
\bibitem{Paper2} \textsc{W.\ Hochst\"attler and M.\ Wilhelmi} \textit{Vertex Shellings of Euclidean Oriented Matroids, preprint}
\bibitem{nelson} \textsc{P.\ Nelson} \textit{Almost all matroids are non-representable, Bull. London Math. Soc. 50 (2018), 245 -- 248.}
\bibitem{SturmfelsZiegler} \textsc{B.\ Sturmfels and G.\ Ziegler} \textit{Extension spaces of oriented matroids. Discrete Comp. Geom. 10, pp 23 - 45, 1993}
 
\end{thebibliography}

\end{document}